\documentclass[10pt]{article}
\usepackage{mathptmx}

\usepackage{color}
\usepackage{footmisc}

\input{epsf.tex}
\usepackage{amsmath,amsthm,amssymb,textcomp}
\usepackage{eucal}
\usepackage{graphicx}
\usepackage{url}
\usepackage[colorlinks,citecolor=blue]{hyperref}

\numberwithin{equation}{section} 

\input xy
\xyoption{all}

\theoremstyle{plain}
\newtheorem*{theo}{Theorem}
\newtheorem{theorem}{Theorem}
\newtheorem{lemma}[theorem]{Lemma}
\newtheorem*{lem}{Lemma}

\newtheorem*{boundedimagetheorem}{Bounded Image Theorem}
\newtheorem*{covering}{Covering Theorem}
\newtheorem*{vol}{Volume Estimate}
\newtheorem*{filling}{Filling Theorem}
\newtheorem*{qf}{Quasifuchsian Density}
\newtheorem*{flex}{Flexible Manifolds}

\newtheorem*{constricting}{Focused Manifolds}
\newtheorem*{pinched}{Pinched Manifolds}

\newtheorem{proposition}[theorem]{Proposition}

\newtheorem*{claim}{Claim}
\newtheorem*{statement}{Statement}

\theoremstyle{definition}
\newtheorem*{question}{Question}

\newtheorem*{problem}{Problem}

\newcommand{\cone}{\mathbf{c}}

\newcommand{\AH}{\mathrm{AH}}
\newcommand{\GF}{\mathrm{GF}}
\newcommand{\QF}{\mathrm{QF}}
\newcommand{\Hom}{\mathrm{Hom}}

\newcommand{\SL}{\mathrm{SL}}
\newcommand{\PSL}{\mathrm{PSL}}

\newcommand{\dee}{\mathrm{d}}
\newcommand{\diam}{\mathrm{diam}}

\newcommand{\co}{\colon\thinspace}

\newcommand{\HH}{\mathbb H}
\newcommand{\R}{\mathbb R}

\newcommand{\C}{\mathbb C}

\newcommand{\F}{\mathfrak F}
\newcommand{\K}{\mathfrak K}
\newcommand{\T}{\mathcal T}
\newcommand{\U}{\mathcal U}
\newcommand{\ELL}{\mathcal L}
\newcommand{\X}{\mathfrak X}

\newcommand{\infinity}{{\rotatebox{90}{\hspace{-.75pt}\large $8$}}}
\newcommand{\smallinfinity}{{\hspace{.5pt} \rotatebox{90}{\footnotesize $8$}}}
\newcommand{\scriptinfinity}{{\rotatebox{90}{\scriptsize $8$}}}

\begin{document}

\title{\textbf{Skinning maps}
}
\author{Richard Peabody Kent IV\thanks{This research was partially conducted during the period the author was employed by the Clay Mathematics Institute as a Liftoff Fellow. The author was also supported by a Donald D. Harrington Dissertation Fellowship and an NSF Postdoctoral Fellowship. }}

\date{June 8, 2009}

\maketitle

\begin{enumerate}
\item[] 
 \small{On this side he fell down from Heaven, and the land which before stood out here made a veil of the sea for fear of him and came to our hemisphere; and perhaps to escape from him that which appears on this side left here the empty space and rushed upwards. }

{\small \quad \quad 
\quad \quad  \quad \quad \quad \quad ---Dante's \textit{Inferno} XXXIV 121--126 (J. D. Sinclair's translation \cite{dante})}
\end{enumerate}





\bigskip


\noindent 
Let $M$ be a compact oriented irreducible atoroidal $3$--manifold with connected incompressible boundary that is not empty nor a torus. 
By W. Thurston's Geometrization Theorem for Haken Manifolds (see \cite{morgansurvey,otalhaken,otalfibered,kapovichbook}), the interior $M^{\circ}$ of $M$ admits a geometrically finite hyperbolic structure with no cusps---G. Perelman has announced a proof of Thurston's Geometrization Conjecture in its entirety \cite{perelman1,perelman2,perelman3}. 
By the deformation theory of L. Ahlfors and L. Bers and stability theorems of A. Marden and D. Sullivan, the space of all such structures $\GF(M)$ may be identified with the Teichm\"uller space $\T(\partial M)$ of $\partial M$ by identifying a manifold with its conformal boundary. 
The cover of $M^\circ$ corresponding to $\partial M$ is a quasifuchsian manifold and we obtain a map
\[
\GF(M) \cong \T(\partial M) \to \T(\partial M) \times \T(\overline{\partial M}) \cong \GF(\partial M \times \R) ,
\]
whose first coordinate function is the identity map on $\T(\partial M)$.
The second coordinate is Thurston's \textbf{skinning map}
\[
\sigma_M \co \T(\partial M) \to \T(\overline{\partial M}),
\]
which reveals a surface obscured by the topology of $M$.

W. Haken proved \cite{haken} that a Haken manifold may be cut progressively along incompressible surfaces until a union of finitely many balls is attained, and Thurston's proof of the Geometrization Theorem proceeds by induction. 
The interior of such a union admits a hyperbolic metric, and this is the first step. 
Thurston reduces each subsequent step to the final one by the so--called Orbifold Trick, turning the part of the boundary one is not gluing into an array of mirrors. 

At the final step you have two manifolds $M$ and $N$, say, each with incompressible boundary homeomorphic to $S$, and a gluing map $\tau \co \partial M \to \partial N$ that produces the manifold under study---there is another case where parts of the boundary of a connected manifold are identified.
Suppose further that $M$ and $N$ are not both interval bundles over a surface---in that case, Thurston's Double Limit Theorem  \cite{thurstongeomII} is used, see \cite{otalfibered}. 
The manifolds $M$ and $N$ have geometrically finite hyperbolic structures given by Riemann surfaces $X$ and $Y$, and the Maskit Combination Theorems \cite{maskit} say that the gluing map produces a hyperbolic structure if ${\tau_{*}}^{\! \! \! -1}(Y) = \sigma_M(X)$ \textit{and} $\tau_*(X) = \sigma_N(Y)$, where $\tau_*$ is the induced map on Teichm\"uller spaces. 
 Such a pair $(X,Y)$ is a \textbf{solution to the gluing problem}.

When the pieces $M$ and $N$ are acylindrical, a solution is provided by  
the following remarkable theorem---when one of the pieces is cylindrical, something else must be done, see Section \ref{bddimage}.
\begin{boundedimagetheorem}[Thurston \cite{thurstonbounded}] If $M$ is acylindrical, the image of $\sigma_M$ has compact closure.\footnote[2]{The statement appears in \cite{thurstonbounded}, Peter Scott's notes (see \cite{scottbangor}) of Thurston's lectures at the conference \textit{Low Dimensional Topology}, in Bangor, in 1979.
The author is grateful to Dick Canary for pointing this out.}
\end{boundedimagetheorem}
\noindent 
Thurston proved that the space $\AH(M)$ of all hyperbolic structures on $M^\circ$ with the algebraic topology is compact \cite{thurstongeomI}, and it can be shown that $\sigma_M$ admits a continuous extension to all of $\AH(M)$---as the proof that $\sigma_M$ admits such an extension has never appeared in print, we supply a proof in Section \ref{bddimage}.
Compactness of $\AH(M)$ was subsequently proven by arboreal methods by J. Morgan and P. Shalen \cite{morganshalen,morganshalen2,morganshalen3}, and
C. McMullen solved the gluing problem via an analytic study of $\sigma_M$ \cite{skinmcmullen}.

\bigskip
\noindent Aside from the conclusion of the Bounded Image Theorem, little has been known about the image of $\sigma_M$. 
Y. Minsky has asked the following question.

\begin{question}[Minsky]
Given a topological description of $M$, can one give a quantitative bound on the diameter of the skinning map?
\end{question}

Notice that it is not obvious that zero does not serve as an upper bound, though
D. Dumas and the author have shown that it will never do:
\begin{theorem}[Dumas--Kent \cite{dumaskent}]\label{nonconstant} Skinning maps are never constant. \qed
\end{theorem}

The map $\GF(M) \to \GF(\partial M \times \R)$ is the restriction of a map
$
\X(M) \to \X(\partial M)
$
at the level of $\SL_2\C$--character varieties. 
We may restrict this map to the irreducible component $\X_0(M)$ containing the Teichm\"uller space $\T(\partial M)$, whose complex dimension is that of $\T(\partial M)$, namely $- \frac{3}{2}\chi(\partial M)$.
If the skinning map $\sigma_M$ were constantly equal to $Y$, then the image of $\X_0(M)$ would contain the Bers slice $\T(\partial M) \times \{ Y \} \subset \QF(\partial M)$. 
This is forbidden by the following theorem.
\begin{theo}[Dumas--Kent \cite{dumaskent}] A Bers slice in $\QF(S)$ is not contained in any subvariety of $\X(S)$ of dimension $- \frac{3}{2}\chi(S)$. \qed
\end{theo}

\bigskip
\noindent 
A sufficiently effective proof of Thurston's fixed point theorem has the potential of transforming his Geometrization Theorem from simply the assertion of the existence of a hyperbolic structure into a statement containing geometric information about that structure, and Minsky's question points to the simplest scenario. 
With this in mind, S. Kerckhoff has put forth the following problem \cite{oberwolfach}.
\begin{problem}[Kerckhoff] Find an effective proof of Thurston's fixed point theorem.
\end{problem}

Our reply to Minsky follows---see Section \ref{volume}.

\begin{vol} Let $M$ be a finite volume hyperbolic $3$--manifold with nonempty closed totally geodesic boundary. 
Then there are positive constants $A$, $B$, and $\epsilon$ depending only on the volume of $M$ such that the image of $\sigma_M$ is $\epsilon$--thick and
\[
B \leq \diam(\sigma_M) \leq A .
\]
\end{vol}

\noindent A set in $\T(S)$ is \textbf{$\epsilon$--thick} if each of its points has injectivity radius at least $\epsilon$.
A theorem of T. J\o rgensen states that the only manner in which a collection of volumes of compact hyperbolic manifolds with totally geodesic boundaries can accumulate is via a sequence of hyperbolic Dehn fillings, see Section \ref{volume}.
The estimates are then obtained as a consequence of the Bounded Image Theorem, Theorem \ref{nonconstant}, and a uniform filling theorem (proven in Section \ref{fillingsection}). 

\begin{filling} Let $M$ be a finite volume hyperbolic manifold with nonempty closed totally geodesic boundary and let $\mathfrak{e} > 0$. 
There is an $\hbar > 0$ such that if the normalized length of each component of a Dehn filling slope $\alpha$ is at least $\hbar$, then
\[
\dee_{\T(\overline{\partial M})}\big(\sigma_M(X), \sigma_{M(\alpha)}(X) \big) < \mathfrak{e}
\]
for all $X$ in $\T(\partial M)$.
\end{filling}
\noindent In other words, as one performs higher and higher Dehn fillings on $M$, the skinning maps of the filled manifolds converge \textit{uniformly} on all of Teichm\"uller space to the skinning map of $M$. 

The reader will notice that we only require the components of $\alpha$ to be long in the cusp cross section of $M$, rather than in all of the geometrically finite manifolds $M_X$ as $X$ ranges over $\T(\partial M)$. 
As it turns out, a ``flat" version of the Bounded Image Theorem tells us that the normalized lengths of a curve in all of the $M_X$ are comparable, and so the normalized length condition is in fact a topological one.
With the Universal Hyperbolic Dehn Filling Theorem of K. Bromberg, C. Hodgson, and S. Kerckhoff, this shows that when $M$ has a single cusp, there is a finite set of filling slopes outside of which hyperbolic Dehn filling may be performed on \textit{all} of the $M_X$.  
There is a similar statement in the presence of a number of cusps, though we warn the reader that it is not always the case that the set of filling slopes exceptional for hyperbolic Dehn filling is finite.  
For example, the Whitehead link complement admits infinitely many Lens space fillings. 
See Section \ref{dehnfilling}.

If one simply desires the skinning maps to be close on a compact set, one may use J. Brock and K. Bromberg's Drilling Theorem \cite{drilling}, and make a geometric limit argument.  
This argument proceeds as follows.  
Given a geometrically finite hyperbolic manifold $M$ with a rank--two cusp and nonempty connected conformal boundary, consider the quasifuchsian group corresponding to the latter.  
The image of its convex core in $M$ misses a neighborhood of the cusp.
By work of Bromberg, hyperbolic Dehn fillings may be performed while fixing the conformal boundary and higher and higher fillings produce bilipschitz maps outside a neighborhood of the cusp with better and better constants, by the Drilling Theorem.  
 Since the core of the quasifuchsian group misses this neighborhood, the generators of this group move less and less in $\PSL_2 \C$.  
 This means that  the quasifuchsian group at $\partial M$ is close in the algebraic topology to the quasifuchsian group at the boundary of the filled manifold. 
 Since the Teichm\"uller metric on the space of quasifuchsian groups induces the same topology, we know that the skinning surfaces are close.
 Notice that this proves that the skinning maps converge.

This argument passes from bilipschitz control to the algebraic topology to the \linebreak Teichm\"uller metric. 
If one attempts the argument as the conformal boundary changes, the intermediate step is an obstacle: as the conformal boundary diverges, one may in fact need better and better quality bilipschitz maps in order to ensure reasonable quality estimates in the Teichm\"uller metric.

To gain uniform control over Teichm\"uller space we must make a more universal argument, which is roughly as follows.  
We have the arrangement above, now with a varying ideal conformal structure $X$.  
Rather than attempting passage through the algebraic topology, we estimate the Teichm\"uller distance between the skinning surfaces directly. 

By the work of C. Epstein, there is a smooth strictly convex surface $\F$ in the end of the quasifuchsian manifold $\QF(X,Y)$ facing $Y$.  
The curvatures of this surface are well behaved and the Hausdorff distance between $\F$ and the convex core is bounded, both independent of $X$.  

There is a universal Margulis tube about the cusp in $M$ that misses the image of $\F$. 
Now, in \cite{drilling}, Brock and Bromberg estimate the strain that the filling process places on $M$ and $\F$
and it follows  that for high enough fillings, the principal curvatures of $\F$ in the filled manifolds are uniformly close to those of $\F$, and so $\F$ is eventually strictly convex in the quasifuchsian covers of these manifolds---if $\hbar$ is large enough, the surface $\F$ lifts to a surface embedded in these quasifuchsian manifolds, see Section \ref{fillingsection}.  

Since $\F$ is convex, there is a normal projection from $\sigma_M(X)$ to $\F$.
Since the image of $\F$ is convex, there is a normal projection from the skinning surface of the filled manifold to the image of $\F$. 
So we obtain a map between skinning surfaces. 
As we do higher and higher fillings, the principal curvatures of the image of $\F$ converge in a controlled way to those of $\F$, and since the derivatives of the projections depend only on these curvatures, they ``cancel in the limit," and we conclude that our map is very close to conformal.  
This yields the desired estimate on the Teichm\"uller distance.

With some additional work, these techniques yield the following theorem due jointly to K. Bromberg and the author---see Section \ref{pinch}.
\begin{pinched}[Bromberg--Kent] Let $S$ be a closed hyperbolic surface. 
For each $\epsilon > 0$ there is a $\delta > 0$ such that if $M$ is an orientable hyperbolic $3$--manifold with totally geodesic boundary $\Sigma$ homeomorphic to $S$ containing a pants decomposition $\mathcal{P}$ each component of which has length less than $\delta$, then the diameter of the skinning map of $M$ is less than $\epsilon$.
\end{pinched}
\noindent In particular, the lower bound in our Volume Estimate tends to zero as the volume grows. 

The proof of the theorem is roughly as follows. 
Let $\Gamma$ be the uniformizing Kleinian group for $M$. 
Since the pants decomposition $\mathcal{P}$ is very short, the Drilling Theorem allows it to be drilled from $\HH^3/\Gamma$.  
This drilled manifold $N$ contains an isometrically embedded copy of the convex core $\mathcal{C}$ of the so--called maximal cusp hyperbolic structure on $M^\circ$ corresponding to $\mathcal{P}$. 
The conformal boundary of the maximal cusp is a union of thrice--punctured spheres, and by a theorem of C. Adams \cite{adams}, these spheres may be taken to be totally geodesic in all of the complete hyperbolic structures on $N$. 
As thrice--punctured spheres have no moduli, the manifold $\mathcal{C}$ isometrically embeds in every hyperbolic structure on $N$.

The meridional Dehn filling slopes of $N$ have large normalized lengths in every hyperbolic structure on $N$, provided $\mathcal{P}$ is short enough, and so Dehn filling may be performed there over all of $\GF(N)$. 
So, by drilling, performing a quasiconformal deformation, and filling, we may convert $\HH^3/\Gamma$ to any point in $\GF(M)$. 

Now a surface $\mathfrak{G}$ with principal curvatures strictly greater than $-1$ is constructed in our copy of $\mathcal{C}$ in $N$.
Unlike the surface used in the proof of the filling theorem, the surface $\mathfrak{G}$ is not convex.  

Pushing $\mathfrak{G}$ back into $\HH^3/\Gamma$ yields a surface whose principal curvatures are still above $-1$, which allows a normal projection to the skinning surface.
As our rigid $\mathcal{C}$ in $N$ is unaffected by any quasiconformal deformation of $N$, so unaffected is our surface $\mathfrak{G}$. 
Filling again still has little affect on $\mathfrak{G}$, and as in the proof of the Filling Theorem, we obtain a map between the skinning surface at $\Sigma$ and the skinning surface at any other point that is very close to being conformal. 
This is the statement that the diameter of $\sigma_M$ is small.

%

\bigskip
\noindent There is a strong form of Minsky's question:
\begin{question}[Minsky] Is the diameter of $\sigma_M$ bounded above by a constant depending only on the topology of $\partial M$?
\end{question}
\noindent A theorem of A. Basmajian \cite{basmajian} implies that a bound on the volume of a manifold with totally geodesic boundary yields a bound on the sum of the genera of the boundary components, and so an affirmative answer implies the upper bound in our Volume Estimate.  
Note that by the previous theorem, and M. Fujii and T. Soma's Density Theorem \cite{fujiisoma}, there is no lower bound depending only on the topology of $\partial M$.

In their solution of Thurston's Ending Lamination Conjecture \cite{ELC,ELC2,ELC3}, Brock, R. Canary, and Minsky show that Minsky's model manifold 
associated to ending laminations $\lambda_-$ and $\lambda_+$ on a surface $S$ is $L$--bilipschitz to any hyperbolic manifold homeomorphic to $S \times \R$ with those ending laminations---this is part of the Bilipschitz Model Theorem. 
It follows from D. Sullivan's Rigidity Theorem \cite{sullivanrigidity} that there is a unique hyperbolic manifold homeomorphic to $S \times \R$ with ending laminations $\lambda_-$ and $\lambda_+$. 
The model is particularly good in this case as the constant $L$ only depends on the topology of $S$.
\begin{figure}
\begin{center}
\includegraphics{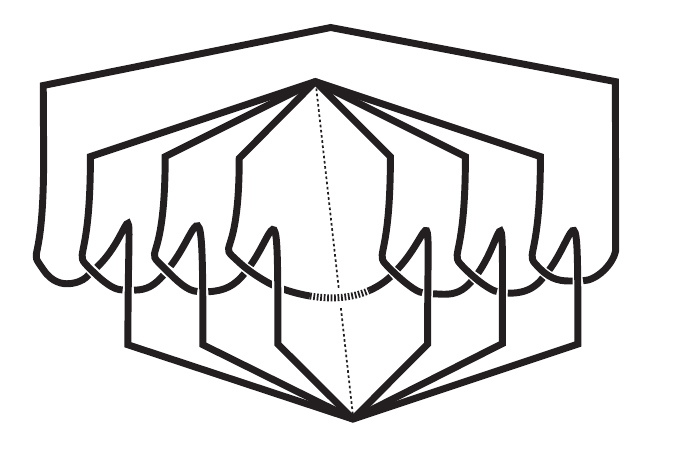}
\end{center}
\caption{Suzuki's Brunnian graph $\mathcal{G}_n$ on $n$ edges.}
\label{graphs}
\end{figure}

In the acylindrical case, the proof of the Model Theorem should proceed as follows.
We have an acylindrical manifold $M$ and a hyperbolic structure on its interior. 
The model is to be built by passing to the cover corresponding to $\partial M$, constructing the model for the end invariants $\lambda$ and $\sigma_M(\lambda)$, cutting off the end facing the Riemann surface $\sigma_M(\lambda)$, and capping off with $M$ equipped with an arbitrary metric.
A uniform bound on the diameter of $\sigma_M$ would be useful in promoting this model to one whose bilipschitz constant depends only on the topology of $\partial M$.

\begin{figure}
\begin{center}
\includegraphics{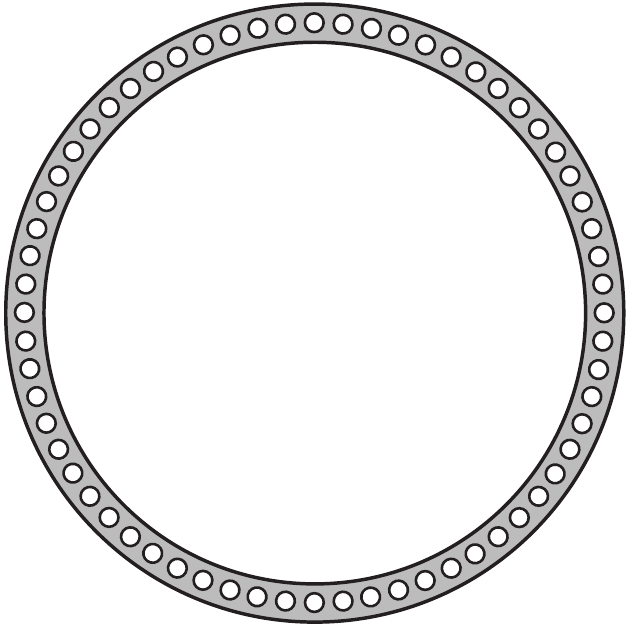}
\end{center}
\caption{The domain of discontinuity of the uniformizing Kleinian group for $S^3 - \mathcal{G}_n$.}
\label{domain}
\end{figure}

As we prove in Section \ref{nouniversal}, any uniformity in an upper bound will necessarily depend on the topology of $\partial M$:

\begin{flex}
There are hyperbolic $3$--manifolds with connected totally geodesic boundary whose skinning maps have diameter as large as you like.
\end{flex}
\noindent The exterior $M_n$ of S. Suzuki's Brunnian graph $\mathcal{G}_n$ in Figure \ref{graphs} \cite{suzuki} admits a hyperbolic structure with totally geodesic boundary, see \cite{paoluzzizimmermann} and \cite{ushijima}. 
As $n$ tends to infinity, the length of the fine arc depicted tends to zero---one may see this directly from the explicit hyperbolic structure given in  \cite{paoluzzizimmermann} and \cite{ushijima}; and indirectly by observing that the $M_n$ are $n$--fold branched covers of $\frac{n}{0}$--orbifold fillings on the complement of a fixed tangle, one strand of which lifts to the fine arc.  
This means that we may normalize the domain of discontinuity to appear as in Figure \ref{domain}, where the disk containing infinity unifomizes the boundary of $M_n$.  
The large disk in the center leads one to suspect that $\sigma_{M_n}$ has large diameter, as a deformation of $\partial M_n$ is carried by the skinning map to a deformation of $\overline{\partial M}_n$ that seems dominated by its effect in the central disk. 
For a particular deformation, this intuition is justified by E. Reich and K. Strebel's First Fundamental Inequality.

In fact, a manifold will have a skinning map of large diameter provided it possesses only shallow collars about its totally geodesic boundary, see Section \ref{nouniversal}; conversely, we have the following theorem, proven in Section \ref{focusing}:

\begin{constricting} Let $S$ be a closed hyperbolic surface.  
Let $\epsilon , \delta, R > 0$. 
There is a constant $d > 0$ such that if $M$ is a hyperbolic $3$--manifold with totally geodesic boundary $\Sigma$ homeomorphic to $S$ of injectivity radius at least $\delta$ possessing a collar of depth $d$, then $\diam \big( \sigma_M\big( B(\Sigma,R) \big) \big) < \epsilon$.
\end{constricting}

A theorem of Fujii and Soma \cite{fujiisoma} says that the surfaces appearing as the totally geodesic boundary of a hyperbolic $3$--manifold without cusps form a dense subset of the Teichm\"uller space---this theorem depends on a theorem of R. Brooks \cite{brooks} that says that the set of surfaces appearing as nonseparating totally geodesic surfaces in closed hyperbolic manifolds is dense.  
It is not difficult to see from the proof that for any $d > 0$, the totally geodesic boundaries of manifolds possessing a collar of depth $d$ about the boundary are also dense, see Section \ref{ubiquity}, and so the manifolds in the theorem exist in great abundance.

To prove the theorem, the skinning surface is compared directly to the totally geodesic boundary as in the proof of the Filling Theorem.  
The necessary control is again guaranteed by work of Brock and Bromberg, namely their Geometric Inflexibility Theorem and the estimates of strain in its proof.

A corollary of the above theorem is the following quasifuchsian version of Fujii and Soma's density theorem, proven in Section \ref{quasifuchsiandensitysection}---again this has the flavor of a theorem of Brooks: the set of quasifuchsian groups induced by nonseparating surfaces in closed hyperbolic manifolds in dense in the space of all such groups, see \cite{brooks}.
\begin{qf} The subset of  $\QF(S)$ consisting of quasifuchsian groups that appear as the only boundary subgroup of a geometrically finite acylindrical hyperbolic $3$--manifold is dense. 
\end{qf}

Another corollary is our Capping Theorem, proven in Section \ref{capsection}, which allows certain properties of skinning maps of manifolds with a number of boundary components to be promoted to those of manifolds with connected boundary.

\bigskip
\noindent
We close the introduction by emphasizing that our manifolds with large diameter skinning maps have \textit{only} very shallow collars about their totally geodesic boundaries, while those with highly contracting maps have very deep collars.  
With Basmajian's theorem \cite{basmajian} that totally geodesic surfaces have collars whose depth depends only on the surface's area, this is perhaps an indication that Minsky's second question has an affirmative answer.

\bigskip

\noindent \textbf{Acknowledgments.}  The author thanks 
Ian Agol, 
Jeff Brock, 
Ken Bromberg, 
Danny Calegari,
Dick Canary, 
Jason DeBlois, 
David Dumas, 
John Holt, 
Chris Leininger, 
Yair Minsky,
Alan Reid, 
Rich Schwartz, 
Juan Souto,
and Genevieve Walsh 
for nice conversation. 
He is especially grateful to Ken Bromberg for his interest, advice, and for allowing the inclusion of Theorem \ref{brombergkent}.
Special thanks are also extended to Jeff Brock and Yair Minsky for conversations that produced the proof of the Bounded Image Theorem given in the final section.

The author also extends his thanks to the referees for their careful readings, comments, and suggestions.

\section{The landscape}

\subsection{Teichm\"uller spaces}

We refer the reader to \cite{ahlfors}, \cite{gardiner}, and \cite{gardlakic} for more on quasiconformal mappings and Teichm\"uller theory.

A homeomorphism $f \co X \to Y$ between Riemann surfaces is \textbf{$K$--quasiconformal} if it has locally integrable distributional derivatives that satisfy
\[
|f_{\bar z}| \leq k |f_z|
\]
where $k <1$ is the number $(K-1)/(K+1)$. 
The \textbf{dilatation} $\mathrm K[f]$ of $f$ is the infimum of all $K$ for which $f$ is $K$--quasiconformal. 
A quasiconformal map has a \textbf{Beltrami coefficient} $\mu_f$ given by the \textbf{Beltrami equation}
\begin{equation}\label{beltrami}
f_{\bar z} = \mu_f f_z 
\end{equation}

Let $\Gamma$ be a torsion free Fuchsian group acting on $\HH^2$---we consider $\HH^2$ as a round disk in $\widehat \C$ of area $\pi$. 
An element $\mu$ of $L^\smallinfinity(\HH^2)$ is a \textbf{Beltrami differential for $\Gamma$} if it satisfies
\[
\mu(g(z)) \frac{\overline{g'(z)}}{g'(z)} = \mu(z)
\]
for all $g$ in $\Gamma$---so that $\mu$ is the lift of a differential $\mu(z) \dee\bar z /\dee z$ on $X_\Gamma = \HH^2/\Gamma$. 
The set of all such differentials is written $L^\smallinfinity(\HH^2, \Gamma)$ or $L^\smallinfinity(X_\Gamma)$. 
Given such a $\mu$ with $\| \mu \|_\smallinfinity < 1$, there is a unique quasiconformal map $f_\mu \co \HH^2 \cup S^1_\smallinfinity \to \HH^2 \cup S^1_\smallinfinity$ fixing $0$, $1$, and $\infinity$ whose Beltrami coefficient is $\mu$. 
Two differentials $\mu$ and $\nu$ are \textbf{Teichm\"uller equivalent} if $f_\mu$ and $f_\nu$ are identical on $S^1_\smallinfinity$---we also say that $f_\mu$ and $f_\nu$ are equivalent in this case. 
The quotient of the open unit ball of $L^\smallinfinity(\HH^2, \Gamma)$ by this equivalence relation is the \textbf{Teichm\"uller space $\T(\Gamma)$ of $\Gamma$}---also referred to as the Teichm\"uller space $\T(X_\Gamma)$ of $X_\Gamma$. 
O. Teichm\"uller's metric on $\T(\Gamma)$ is given by
\[
\dee_{\T(\Gamma)}\big([\mu], [\nu]\big) = \frac{1}{2} \inf \, \log \mathrm{K}[f_{\widetilde \mu} \circ f_{\widetilde \nu}^{-1}]
\]
where the infimum is over all $\widetilde \mu$ in $[\mu]$ and $\widetilde \nu$ in $[\nu]$.

A quasiconformal $f$ is \textbf{extremal} if the $L^\smallinfinity$--norm of its Beltrami coefficient is less than or equal to that of any map equivalent to $f$. 
Every quasiconformal map is equivalent to an extremal one. 
Given a Beltrami differential $\mu$ in the unit ball of $L^\smallinfinity(X_\Gamma)$, we let $\mathrm K[\mu]_{X_\Gamma}$ denote the dilatation of an extremal quasiconformal map equivalent to $f_\mu$, called its \textbf{extremal dilatation}, and we drop the subscript $X_\Gamma$ when no confusion can arise.

We write $\mathcal{Q}(X_\Gamma)$ for the space of integrable holomorphic quadratic differentials on $X_\Gamma$, equipped with the norm
\[
\| \, \cdot \, \| = \int_{X_\Gamma} | \, \cdot \, | 
\]
On $\HH^2$, these are the integrable holomorphic functions $\varphi$ satisfying
\[
\varphi(g(z)) g'(z)^2 = \varphi(z)
\]
for all $g$ in $\Gamma$ and all $z$ in $\HH^2$---these are the \textbf{holomorphic cusp forms for $\Gamma$}.

\subsubsection{The Fundamental Inequality}

The inequality here due to Reich and Strebel estimates the extremal dilatation of a Beltrami differential from below.
We refer the reader to Section 4.9 of \cite{gardlakic} for a proof.

\begin{theorem}[First Fundamental Inequality \cite{gardlakic}] 
Let $f \co X \to Y$ be a quasiconformal map with Beltrami coefficient $\mu$.
Then the extremal dilatation satisfies 
\begin{equation}\label{firstineq}
\frac{1}{\mathrm K[\mu]_X}  
\leq   \iint_{X} |\varphi|
\frac{\ \big| 1 - \mu \frac{\varphi}{| \varphi |} \big|^2}{1-|\mu|^2}
\end{equation}
for all unit norm holomorphic quadratic differentials $\varphi$ on $X$. \qed
\end{theorem}

\noindent Note that if $\mu$ is a Teichm\"uller differential $k | \varphi |/ \varphi$, we have equality in \eqref{firstineq}.

\subsection{$3$--manifolds}

Let $M$ be a smooth compact $3$--manifold with boundary.
The manifold $M$ is \textbf{irreducible} if every smoothly embedded $2$--sphere bounds a ball. 
It is \textbf{atoroidal} if every $\pi_1$--injective map $S^1 \times S^1 \to M$ is homotopic into $\partial M$, and \textbf{acylindrical} if every homotopically essential map of pairs $(A, \partial A) \to (M,\partial M)$ is homotopic as a map of pairs into $\partial M$ whenever $A$ is an annulus $S^1 \times I$. 
A noncompact $3$--manifold is \textbf{tame} if it is homeomorphic to the interior of a compact manifold.

\subsection{Kleinian groups}

The group of orientation preserving isometries of hyperbolic $3$--space $\HH^3$ is $\PSL_2\C$, and a \textbf{Kleinian group} is a discrete subgroup of this group. 
We will assume throughout that our Kleinian groups are \textbf{nonelementary}, meaning that they are not virtually abelian.

The \textbf{limit set} $\Lambda_\Gamma$ of a Kleinian group $\Gamma$ is the unique minimal nonempty closed $\Gamma$--invariant subset of $\widehat \C$.  
Its complement, $\Omega_\Gamma = \widehat \C - \Lambda_\Gamma$ is the \textbf{domain of discontinuity} of $\Gamma$.

We let $M_\Gamma = \HH^3/\Gamma$ and let $\dot M_\Gamma = (\HH^3 \cup \Omega_\Gamma)/\Gamma$ denote the \textbf{Kleinian manifold}.
The Riemann surface $\Omega_\Gamma/\Gamma$ is the \textbf{conformal boundary} of $M_\Gamma$.
By the Ahlfors Finiteness Theorem \cite{ahlforsfiniteness,ahlforsfinitenesscorrection,greenberg}, see also \cite{kapovichfiniteness}, the conformal boundary of a finitely generated Kleinian group has finite type.

The quotient $\mathcal{C}_\Gamma = \mathfrak{H}_\Gamma/\Gamma$ of the convex hull $\mathfrak{H}_\Gamma$ of $\Lambda_\Gamma$ is the \textbf{convex core} of $M_\Gamma$. 
We sometimes write $\mathcal{C}(M_\Gamma)$ for $\mathcal{C}_\Gamma$. 
We say that $\Gamma$ and $M_\Gamma$ are \textbf{geometrically finite} if the $1$--neighborhood of $\mathcal{C}_\Gamma$ has finite volume.

Given a hyperbolic manifold $M$, we write $\Gamma_M$ for the uniformizing Kleinian group.

We let $\epsilon_3$ be the $3$--dimensional Margulis constant.

\subsubsection{The algebraic topology}

Let $M$ be a compact oriented $3$--manifold with \textit{incompressible boundary} and let $P$ be a $\pi_1$--injective $2$--dimensional submanifold of $\partial M$ of Euler characteristic zero containing all of the tori in $\partial M$.
The \textbf{space of hyperbolic structures} $\mathrm{H}(M, P)$ is the space of discrete faithful representations $\rho \co \pi_1(M) \to \PSL_2 \C$ up to conjugacy such that for each component $P_0$ of $P$, the image $\rho\big(\pi_1(P_0)\big)$ is purely parabolic. 
Such representations are then holonomy representations of complete hyperbolic manifolds homotopy equivalent to $M$ and so we may think of $\mathrm{H}(M,P)$ as a space of hyperbolic manifolds.
We let $\AH(M,P)$ denote this set equipped with the topology induced by the inclusion
\[
\AH(M, P) \subset \Hom(\pi_1(M), \PSL_2 \C)\big/\PSL_2 \C \, ,
\]
called the \textbf{algebraic topology}---this is the same topology that $\AH(M,P)$ inherits from the inclusion into the $\SL_2\C$--character variety $\X(M,P)$, which is a particular birational representative of the quotient by $\SL_2\C$ of the space of representations with parabolics at $P$ in the sense of geometric invariant theory (we refer the reader to \cite{cullershalen} and \cite{shalen} for detailed discussions of the $\SL_2\C$--character variety). 
When $P$ is the union of all of the tori in $\partial M$, we write $\AH(M) = \AH(M,P)$. 
When given a geometrically finite hyperbolic $3$--manifold $N$, possibly with geodesic boundary, we often write $\AH(N)$ for the space $\AH(M, P)$ where $M$ is a compact manifold whose interior is homeomorphic to that of $N$ and $P$ is the submanifold of $\partial M$ corresponding to the cusps of $N$.
A hyperbolic structure in $\AH(M,P)$ has \textbf{totally geodesic boundary} if  the subgroups of $\pi_1(M)$ corresponding to $\partial M - P$ are all Fuchsian, so that the representation is the holonomy representation of a hyperbolic metric on $M - P$ with totally geodesic boundary.

Combining work of V. Chuckrow \cite{chuckrow} and J\o rgensen \cite{jorgensen}, it follows that $\AH(M,P)$ is closed in $\X(M,P)$, and, as mentioned in the introduction, the key to the Bounded Image Theorem is the following compactness theorem of Thurston, see also \cite{morganshalen,morganshalen2,morganshalen3}.
\begin{theorem}[Thurston \cite{thurstongeomI,thurstongeomIII}]\label{AHcompact} If $(M,P)$ admits a hyperbolic structure with totally geodesic boundary, then $\AH(M,P)$ is compact. \qed
\end{theorem}

By a theorem of Marden \cite{marden} and Sullivan \cite{sullivanstability}, one component of the interior of $\AH(M,P)$ is precisely the set of geometrically finite hyperbolic structures on $M^\circ$ compatible with the orientation on $M$ with no parabolics other than those arising from $P$. 
We let $\GF(M,P)$ denote this component. 
Identifying a manifold with its conformal boundary, the Measurable Riemann Mapping Theorem of Ahlfors and Bers \cite{ahlforsbers} then implies the following theorem, see \cite{bersparameter} or \cite{canarymccullough}.

\begin{theorem}[Ahlfors--Bers--Marden--Sullivan]\label{geomfinite} Quasiconformal conjugation induces a biholomorphic isomorphism 
$\GF(M,P) \cong \T(\partial M - P)$. \qed
\end{theorem}
\noindent Given a point $X$ in $\T(\partial M - P)$, we write $M_X$ for the manifold $M^\circ$ equipped with the hyperbolic structure corresponding to $X$.

Let $S$ be a closed surface. 
When $M$ is homeomorphic to $S \times [0,1]$, we write $\QF(S) = \GF(M)$. 
In this case, Theorem \ref{geomfinite} is Bers' Simultaneous Uniformization Theorem:

\begin{theorem}[Simultaneous Uniformization \cite{berssimultaneous}] $\QF(S) \cong \T(S) \times \T(\overline{S})$. \qed
\end{theorem}

If $X$ and $Y$ are points in $\T(S)$ and $\T(\overline S)$ we write $\QF(X,Y)$ for the hyperbolic manifold in $\QF(S)$ corresponding to the pair $(X,Y)$.  
If $A$ and $B$ are sets in $\T(S)$ and $\T(\overline S)$, we let $\QF(A,B)$ denote the set of all manifolds $\QF(X,Y)$ with $X$ in $A$ and $Y$ in $B$.

\subsubsection{The geometric topology}

A sequence $\Gamma_n$ of Kleinian groups converges in the \textbf{Chabauty topology} to a group $\Gamma < \PSL_2\C$ if the following conditions hold
\begin{enumerate} 
\item For each $\gamma$ in $\Gamma$, there are elements $\gamma_n$ in $\Gamma_n$ such that $\gamma = \lim \gamma_n$. 
\item If $\Gamma_{n_k}$ is a subsequence,  and $\gamma_{n_k}$ in $\Gamma_{n_k}$ are such that $\lim \gamma_{n_k} = \gamma$, then $\gamma$ lies in $\Gamma$.
\end{enumerate}

Geometrically, this is formulated as follows.
A \textbf{framed hyperbolic manifold} is simply a hyperbolic manifold equipped with an orthonormal frame at a basepoint.
A sequence of framed hyperbolic manifolds $M_n$ \textbf{converge geometrically} to a framed hyperbolic manifold $M$ if for each smooth compact submanifold $K$ of $M$ containing the basepoint, there is a smooth frame preserving map $\varphi_n \co K \to M_n$, and the $\varphi_n$ converge to an isometry in the $C^{\, \smallinfinity} \! \!$--topology---meaning that the lifted maps $\widetilde \varphi_n \co \widetilde K \to \HH^3$, normalized so that all frames are the standard frame at the origin, converge in the topology of $C^{\, \smallinfinity} \! \!$--convergence on compact sets.

These two notions are equivalent, see Chapter 2 of \cite{renorm} and Chapter E of \cite{bp}.

\subsubsection{Pleated surfaces}

A map $f \co \Sigma \to M$ from a hyperbolic surface $\Sigma$ to a hyperbolic $3$--manifold $M$ is a \textbf{pleated surface} provided it is a path isometry (meaning that it sends rectifiable arcs to rectifiable arcs of the same length) and for each $x$ in $\Sigma$, there is a geodesic segment $\gamma$ through $x$ that $f$ carries to a geodesic segment in $M$.  

Confusing a pleated surface with its image, the boundary $\partial \mathcal{C}_\Gamma$ of the convex core of a hyperbolic manifold $M_\Gamma$ is a collection of pleated surfaces. 

A \textbf{framed pleated surface} is a pleated surface $f \co \Sigma \to M$ together with choices of orthonormal frames at points $p$ and $q$ in $\Sigma$ and $M$, respectively, such that $f$ carries the frame at $p$ into the frame at $q$.

For $A, \epsilon > 0$, let $\mathcal{PSF}(A, \epsilon)$ be the set of all framed pleated surfaces $f \co \Sigma \to M$ into hyperbolic $3$--manifolds $M$ such that $\Sigma$ has area less than $A$ and such that the injectivity radii of $\Sigma$ and $M$ are bounded below by $\epsilon$ at the base frame.
Using the Chabauty topology on the uniformizing Fuchsian and Kleinian groups for the $\Sigma$ and $M$, and the compact--open topology for the maps, we obtain a topology on $\mathcal{PSF}(A, \epsilon)$.

We will need the following theorem of Thurston \cite{thurstongeomI}, see Theorem 5.2.2 of \cite{notesonnotes}.

\begin{theorem}[Pleated Surfaces Compact]\label{pleatedcompact} The space $\mathcal{PSF}(A, \epsilon)$ is compact.  \qed 
\end{theorem}

\bigskip \noindent
We will also need the following standard lemma.
\begin{lemma}\label{brooksmatelskilemma}
Let  $\epsilon_3 \geq \epsilon > 0$ and let $f \co \Sigma \to M$ be a pleated surface.
There is a $\delta_\epsilon < \epsilon$ depending only on $\epsilon$ and the topological type of $\Sigma$ such that if  $f(\Sigma)$ intersects a component  $\mathbf{T}^{\delta_\epsilon}$ of the $\delta_\epsilon$--thin part of $M$, then there is an essential simple closed curve $\gamma$ in $\Sigma$ of length less than $\epsilon$ such that $f(\gamma) \subset \mathbf{T}^{\epsilon}$, where $\mathbf{T}^{\epsilon}$ is the component of the $\epsilon$--thin part containing $\mathbf{T}^{\delta_\epsilon}$.
\end{lemma}
\begin{proof}
A theorem of Brooks and Matelski \cite{brooksmatelski} says that the distance between $\partial \mathbf{T}^\epsilon$ and $\partial \mathbf{T}^\delta$ tends to infinity as $\delta$ tends to zero. 
So, if $f(\Sigma)$ intersects $\mathbf{T}^\delta$ for some small $\delta$, then, since $f$ is $1$--Lipschitz, the set $f^{-1}(\mathbf{T}^{\epsilon})$ either contains an essential simple closed curve of length less than $\epsilon$ or a disk of enormous area.  
The Gauss--Bonnet Theorem completes the proof.
\end{proof}

\subsubsection{Deformations}

A diffeomorphism $f \co (M,g) \to (N,h)$ between Riemannian manifolds is \textbf{$L$--bilipschitz} if there is an $L \geq 1$ such that for every $x$ in $M$,
\[
\frac{1}{L} \, \| v \|_g \leq \| \, D f _x v \, \|_h \leq L \, \| v \|_g
\]
for all $v$ in $T_x M$.
We will make use of the following well--known theorem---see Theorem 2.5 and Corollary B.23 of \cite{renorm}.
\begin{theorem}\label{quasiiso} Let $\Gamma$ be a Kleinian group.  

A $\Gamma$--equivariant $L$--bilipschitz map $\HH^3 \to \HH^3$ extends continuously to a $\Gamma$--equivar-iant $K(L)$--quasiconformal map $\widehat \C \to \widehat \C$. 
Moreover, $K(L)$ tends to $1$ as $L$ does.

A $\Gamma$--equivariant $K$--quasiconformal map $\widehat \C \to \widehat \C$ extends continuously to a $\Gamma$--equivariant $K^{3/2}$--bilipschitz map $\HH^3 \to \HH^3$. \qed
\end{theorem}

\noindent And we will once mention the following proposition---see Proposition 2.16 in \cite{renorm}.

\begin{proposition}\label{movinghulls} Let $\Psi\co \HH^3 \to \HH^3$ be an $L$--bilipschitz map and let $\Lambda$ be a closed set in $\widehat \C$ with convex hull $\mathfrak{H}_\Lambda$.  
Then there is a constant $h(L)$ such that $\Psi(\mathfrak{H}_\Lambda)$ lies in the $h(L)$--neighborhood of $\mathfrak{H}_{\Psi(\Lambda)}$. \qed \end{proposition}

\subsection{Skinning}

Let $M$ be a hyperbolic manifold with connected totally geodesic boundary $\Sigma$.  
There is a map induced by inclusion
\[
\iota \co \T(\Sigma) \cong \GF(M) \longrightarrow \QF(\Sigma) \cong \T(\Sigma) \times \T(\overline{\Sigma}) 
\]
whose first coordinate function is the identity on $\T(\Sigma)$.
Composition with projection onto the second factor yields a map
\[
\sigma_M \co \T(\Sigma) \to \T(\overline{\Sigma}),
\]
called the \textbf{skinning map} associated to $M$.
The fact that the target of the map $\iota$ is indeed $\QF(\Sigma)$ follows from the Annulus Theorem \cite{waldhausenAnnulus,cannonfeustel} and the fact, first observed by Thurston, that a geometrically finite Kleinian group with nonempty domain of discontinuity has the property that all of its finitely generated subgroups are geometrically finite,  which itself follows from the Ahlfors Finiteness Theorem, see \cite{morgansurvey}.
As the map $\iota \co \GF(M) \to \QF(\Sigma)$ is the restriction of the regular map $\X(M) \to \X(\Sigma)$ between $\SL_2\C$--character varieties induced by inclusion, and the identifications $\T(\Sigma) \cong \GF(M)$ and $\QF(\Sigma) \cong  \T(\Sigma) \times \T(\overline{\Sigma})$ are biholomorphic isomorphisms, the skinning map is holomorphic.

We henceforth denote the open unit disk in $\widehat \C$ by $\ELL$, and let $\U$ denote the interior of its complement.

After conjugating, we assume that the boundary $\Sigma$ is uniformized by a Fuchsian subgroup $\Gamma_\Sigma \leq \Gamma_M$  acting on $\U$.

The skinning map is easily described at the level of equivariant Beltrami differentials by sending the Teichm\"uller class of a differential $\mu$ in $L^\smallinfinity(\U, \Gamma_\Sigma)$ to the class of the differential in $L^\smallinfinity(\ELL, \Gamma_{\overline{\Sigma}})$ obtained by transporting $\mu$ to a differential on $\widehat \C$ via the functional equation
\begin{equation}\label{functional}
\mu\big(g(w)\big) \frac{\overline{g'(w)}}{g'(w)} = \mu(w) \quad \mathrm{when}\ g \in \Gamma_M
\end{equation}
and restricting to $\ELL$. 

Suppose that $\Sigma = \Sigma_0 \sqcup \Sigma_1 \sqcup \ldots \sqcup \Sigma_m$, and that $\Sigma_0$ is connected and uniformized by $\Gamma_{\Sigma_0}$ in the disk $\U$.  
We define the \textbf{skinning map $\sigma_{M,\Sigma_0}$ of $M$ relative to $\Sigma_1 \sqcup \ldots \sqcup \Sigma_m$} by transporting a differential in $L^\smallinfinity(\U, \Gamma_\Sigma)$ to $\Gamma_M \U$ using \eqref{functional},  declaring that the differential be identically zero on $\widehat \C - \Gamma_M \U$, and restricting to $\ELL$.

Given a skinning map $\sigma$, we will, to simplify the notation, write $\sigma(\mu)$ for the differential in $L^\smallinfinity(\ELL, \Gamma_{\overline{\Sigma}})$ constructed above, and $\sigma([\mu]) = [\sigma(\mu)]$ for its Teichm\"uller equivalence class---although the latter is the actual image of $[\mu]$ under the skinning map, it is useful to work with the representative $\sigma(\mu)$.

\subsection{Dehn filling for families}\label{dehnfilling}

Given a  $3$--manifold $M$, a \textbf{slope} is the isotopy class of an essential embedded $1$--manifold on the union of tori in $\partial M$ with at most one component on each torus. 
Given a slope $\alpha$, we may perform \textbf{Dehn filling} along $\alpha$ to obtain a manifold $M(\alpha)$---if $\alpha$ is a union of circles $\alpha_j$, then $M(\alpha)$ is obtained by gluing solid tori $T_j$ to $\partial M$ so that $\alpha_j$ bounds a disk in $T_j$.

If $M$ is a geometrically finite hyperbolic manifold, \textbf{slope} will mean an isotopy class of an essential embedded $1$--manifold on the boundary of the complement of a horospherical neighborhood of the cusps, and Dehn fillings will refer to fillings of this complement.  
If the $\alpha$--filling of $M$ admits a complete hyperbolic metric, we reserve $M(\alpha)$ to denote the hyperbolic structure \textit{with the same conformal boundary} as $M$, and a filling is \textbf{hyperbolic} if the cores of the filling tori are geodesics---in J\o rgensen's theorem we alter this convention, see Section \ref{volume}.

In the finite volume setting, Thurston's Hyperbolic Dehn Filling Theorem  provides a set of filling slopes outside of which any filling may be taken to be hyperbolic, see \cite{fillingtheorem} for instance.  
Hodgson and  Kerckhoff have shown \cite{HKuniversal} that as long as the normalized length of each component of a slope $\alpha$ in the horospherical cross sections of the cusps is at least $7.515 \sqrt 2$, then hyperbolic Dehn filling may be performed at $\alpha$---the \textbf{normalized length} of a curve $\gamma$ in a flat torus is its length divided by the positive square root of the area of the torus, the normalized length of an isotopy class the infimum of the normalized lengths of its representatives.

Given a geometrically finite hyperbolic manifold $M$, a generalization of the Hyperbolic Dehn Filling theorem---first used by F. Bonahon and J.-P. Otal \cite{bonahonotal}---applies to show that, at each cusp, there is a finite set of slopes that are exceptional for hyperbolic filling \cite{comar}.  
It is not \textit{a priori} obvious that this set may be taken to be independent of the geometrically finite hyperbolic structure on $M$.\footnote[3]{A slope $\alpha$ that is nonexceptional for hyperbolic Dehn filling on $M$ may be exceptional for hyperbolic Dehn filling on $M_X$.
To see this, note that while the manifold $M(\alpha)$ does admit a hyperbolic structure $M(\alpha)_X$ with conformal boundary $X$, thanks to Thurston's Geometrization Theorem for Haken Manifolds, it is not clear that the core of the filling torus may be taken to be geodesic in $M(\alpha)_X$.
So $\alpha$ may be exceptional for hyperbolic Dehn filling on $M_X$.}  
Thanks to the work of Bromberg, Hodgson and Kerckhoff's Universality Theorem extends (with a different constant) to the geometrically finite setting, which gives to each cusp a uniform finite set of exceptional slopes for all hyperbolic structures in $\GF(M)$, as we now discuss.

\bigskip
\noindent Let $M$ be a compact orientable irreducible acylindrical $3$--manifold. 
Let $T$ be the union of the tori in $\partial M$.  
So, $\GF(M) \cong \T(\partial M - T)$. 
Let $\T(T)$ be the Teichm\"uller space of $T$---the space of marked flat metrics on $T$ such that each component of $T$ has area one. 
There is a map
\[
\tau_{\, M} \co \GF(M) \to \T(T)
\]
that assigns the shapes of the cusp cross sections to a geometrically finite hyperbolic structure. 
We have the following ``flat" version of the Bounded Image Theorem.

\begin{theorem}\label{flatbounded} The image of $\tau_{\, M}$ has compact closure. 
\end{theorem}
\begin{proof}  The space $\AH(M)$ is compact, and  $\tau_M$ extends continuously to the closure $\overline{\GF(M)}$ of $\GF(M)$ in $\AH(M)$.
\end{proof}
\noindent While continuity in the proof is apparent \textit{here}, it is not in the proof of the Bounded Image Theorem, see Section \ref{bddimage}.

The following theorem is a version of Hodgson and Kerckhoff's Universal Hyperbolic Dehn Filling Theorem, extended via work of Bromberg, that allows the manifolds to have infinite volume, and is stated in \cite{brombergnotlocally}.
\begin{theorem}[Bromberg--Hodgson--Kerckhoff \cite{HKuniversal,brombergJAMS,brombergrigidity}]\label{universal}
There is a universal constant $h$ such that the following holds.
If $M$ is as above and the normalized length of each component of a filling slope $\alpha$ is at least $h$ on $T$, then hyperbolic Dehn filling may be performed at $\alpha$. 
Moreover, for each $\ell$, there is an $h_\ell$ such that if these normalized lengths are at least $h_\ell$, then the lengths of the cores of the filling tori are at most $\ell$. \qed
\end{theorem}

\noindent The hyperbolic Dehn filling is achieved by a cone deformation that pushes the cone angle from zero up to $2\pi$, see Section \ref{drillsection}. 

Theorem \ref{universal} has the following topological corollary.

\begin{theorem}[Universal Exceptions]\label{exceptions} Let $M$ be as above and let $\ell > 0$. 
Let $\mathbf{P}_1, \ldots , \mathbf{P}_n$ be the rank--two cusps of $M$, and let $\mathfrak{S}_i$ be the set of slopes on $\partial \mathbf{P_i}$. 
Then there are finite sets $\mathcal{E}_i \subset \mathfrak{S}_i$ such that for all $M_X$ in $\GF(M)$ and all 
\[
\alpha \in \{ 
(\alpha_1, \ldots, \alpha_n) \ | \ \alpha_i \notin \mathcal{E}_i \mathrm{\ for\ all\ } i \,
\}, 
\]
one may perform hyperbolic Dehn filling on $M_X$ at $\alpha$. 
Moreover,
the lengths of the cores of the filling tori are less than $\ell$ in $M_X(\alpha)$ for all $M_X$ in $\GF(M)$. 
\end{theorem}
\begin{proof} 
The theorem is simply a corollary of Theorem \ref{universal}, Theorem \ref{flatbounded}, and the fact that the length spectrum of a flat torus is discrete with finite multiplicities.
\end{proof}

\subsection{Strain fields}\label{strain}

Deformations of hyperbolic manifolds such as cone deformations and quasiconformal conjugacies induce bilipschitz maps between manifolds. 
In these cases, the pointwise bilipschitz constant may be controlled by estimating norms of certain strain fields. 
These estimates often provide control of other geometric information, such as various curvatures of embedded surfaces, which we will soon desire and so now discuss.

The material in the next few sections is drawn mostly from \cite{drilling}, \cite{bbinflex}, and $\cite{renorm}$, where the reader will find more detailed discussions---see also \cite{HKsurvey}.

In the following, $g_t$ will denote a smooth family of hyperbolic metrics on a $3$--manifold $M$. 
There are vector valued $1$--forms $\eta_t$ defined by
\[
\frac{\dee g_t(x,y)}{\dee t} = 2 g_t\big(x, \eta_t(y)\big)
\]
and, choosing an orthonormal frame field $\{ e_1,  e_2, e_3 \}$ for the $g_t$ metric, we define the norm
\begin{equation}\label{pointwisenorm}
\| \eta_t \|^2 = \sum_{i,j=1}^3 g_t\big(\eta_t(e_i), \eta_t(e_j) \big).
\end{equation}

We consider metrics $g_t$ that are the pullback of a fixed metric $g$ under the flow $\varphi_t$ of a time--dependent vector field $v_t$. 
Letting $\nabla^t$ be the Levi--Civita connection for $g_t$, we have
\[
\eta_t = \mathrm{sym}\, \nabla^t v_t
\]
when $g_t = \varphi_t^*g$. 
The $v_t$ will be \textbf{harmonic}, meaning that they are divergence free and $\mathrm{curl\, curl} \, v_t = - v_t$ (where the $\mathrm{curl}$ is half of the usual one). 
On vector--valued $k$--forms we define the operator
\[
D_t = \sum_{i=1}^3 \omega^i \wedge \nabla^t_{e_i}
\]
where $\{ \omega^1, \omega^2, \omega^3 \}$ is a coframe field dual to $\{ e_1, e_2, e_3 \}$. 
The formal adjoint of this operator is then
\[
D_t^* = \sum_{i=1}^3 \iota(e_i) \nabla^t_{e_i}
\]
where $\iota$ is contraction.
In \cite{HK} it is shown that $D_t^* \eta_t = 0$.

\begin{proposition}[Brock--Bromberg, Proposition 6.4 of \cite{drilling}]\label{strainofcurvature} Let $\gamma(s)$ be a smooth curve in $M$ and let $C(t)$ be the geodesic curvature of $\gamma$ at $\gamma(0) = p$ in the $g_t$--metric.  
For each $\epsilon > 0$ there exists a $\delta > 0$ depending only on $\epsilon$, $a$, and $C(0)$ such that $| \, C(a) - C(0) \, | < \epsilon$ if $\| \eta_t(p) \| \leq \delta$, $\| *D_t \eta_t(p) \| \leq \delta$, and $D_t^*\eta_t = 0$ for all $t \in [0,a]$. \qed
\end{proposition}

The pointwise norm of $*D_t \eta_t$ is defined as in \eqref{pointwisenorm}.

Given a vector field $v$, the symmetric traceless part of $\nabla v$ is the \textbf{strain} of $v$, which measures the conformal distortion of the metric pulled back by the flow of $v$.  
When $v$ is divergence--free  $\nabla v$ is traceless, and so $\eta = \mathrm{sym}\, \nabla v$ is a strain field. 
We say that $\eta$ is \textbf{harmonic} if $v$ is.

If $\eta$ is a harmonic strain field, so is $*D\eta$, by Proposition 2.6 of \cite{HK}.

We will need the following mean--value inequality for harmonic strain fields, see Theorem 9.9 of \cite{brombergJAMS}.

\begin{theorem}[Hodgson--Kerckhoff]\label{meanvalue} Let $\eta$ be a harmonic strain field on a ball $B_R$ of radius $R$ centered at $p$.  
Then
\[
\| \eta(p) \| \leq \frac{3 \sqrt{2 \, \mathrm{vol}(B_R)} }{4 \pi f(R)} 
\sqrt{\int_{B_R} \| \eta \|^2 \dee V\ }  , 
\]
where
\[
f(R) = \cosh(R) \sin\big(\sqrt{2} R\big) - \sqrt{2} \sinh(R) \cos\big(\sqrt{2} R\big)
\]
and $R < \pi \big/ \sqrt{2}$. \qed
\end{theorem}

\subsection{Cone manifolds}

Let $M$ be a compact manifold and let $\cone \subset M$ be a compact $1$--manifold.  
A \textbf{hyperbolic cone metric} on $(M,\cone)$ is a hyperbolic metric on the interior of $M - \cone$ whose metric completion is a metric on the interior of $M$ that in a neighborhood of a component of $\cone$ has the form
\[
\dee s^2 = \dee r^2 + \sinh^2 r \, \dee \theta^2 + \cosh^2 r \, \dee z^2
\]
where $\theta$ is measured modulo a \textbf{cone angle} $a$. 

The hyperbolic metric on the universal cover of the interior of $M - \cone$ has a developing map to $\HH^3$, and we say that a family $g_t$ of hyperbolic cone metrics is smooth if its family of developing maps is smooth.

Our hyperbolic cone manifolds will be geometrically finite, meaning that the cone metric extends to a conformal structure on the portion of $\partial M$ that contains no tori.

\subsection{Drilling}\label{drillsection}

The following theorem of Bromberg allows one to deform certain cone manifolds by pushing the cone angle all the way to zero while fixing the conformal boundary.

\begin{theorem}[Bromberg, Theorem 1.2 of  \cite{brombergJAMS}]\label{brombergexistence} Given $a>0$ there exists an $\ell > 0$ such that the following holds.
 Let $M_a$ be a geometrically finite hyperbolic cone manifold with no rank--one cusps, singular locus $\cone$ and cone angle $a$ at $\cone$.  
If the tube radius $R$ about each component of $\cone$ is at least $\mathrm{arcsinh}\, \sqrt 2$ and the total length $\ell_{M_a}(\cone)$ of $\cone$ in $M_a$ satisfies $\ell_{M_a}(\cone) < \ell$, then there is a $1$--parameter family $M_t$ of cone manifolds  with fixed conformal boundary and cone angle $t \in [0,a]$ at $\cone$. \qed
\end{theorem}

\noindent (In our work here, the tube radius condition is easily established, as we always begin or end with cone angle $2 \pi$, and a theorem of Brooks and J. P. Matelski \cite{brooksmatelski} tells us that the radius of a tube is large when its core is short.) Together with the Hodge Theorem of Hodgson and Kerckhoff \cite{HK} and Bromberg's generalization of it to geometrically finite cone manifolds \cite{brombergJAMS}, this yields the following theorem.

\begin{theorem}[Bromberg--Hodgson--Kerckhoff, Corollary 6.7 of \cite{drilling}] Let $M_t$ be the $1$--parameter family in Theorem \ref{brombergexistence}. 
There exists a $1$--parameter family of cone metrics $g_t$ on $M$ such that $M_t = (M, g_t)$ and $\eta_t$ is a harmonic strain field outside a small tubular neighborhood of the singular locus and the rank--two cusps. \qed
\end{theorem} 

Brock and Bromberg establish an $L^2$--bound on the norm of this strain field:

\begin{theorem}[Brock--Bromberg, Theorem 6.8 of \cite{drilling}]\label{L2bounddrilling} Given $\epsilon > 0$, there are $\ell, K > 0$ such that if $\ell_{M_t}(\cone) \leq \ell$, then 
\[
\int_{M_t - \mathbf{T}_t^\epsilon(\cone)} \| \eta_t \|^2 + \| * D_t\eta_t \|^2 
\ \leq \ K^2 \ell_{M_t}(\cone)^2 
\] 
where $\mathbf{T}_t^\epsilon(\cone)$ is the $\epsilon$--Margulis tube about $\cone$ in $M_t$. \qed
\end{theorem}

\noindent This allows them to prove the following theorem.

\begin{theorem}[Brock--Bromberg, Theorem 6.2 of  \cite{drilling}]\label{drilling} For each $\epsilon_3 \geq \epsilon > 0$ and $L > 1$, there is an $\ell > 0$ such that the following holds.
If $M$ is a geometrically finite hyperbolic $3$--manifold and $\cone$ is a geodesic in $M$ with length $\ell_M(c) < \ell$, there is an $L$--bilipschitz diffeomorphism of pairs
\[
h \co \big(M - \mathbf{T}^\epsilon (\cone), \ \partial \mathbf{T}^\epsilon (\cone) \big) 
\longrightarrow \big(M_0 - \mathbf{P}^\epsilon (\cone), \, \partial \mathbf{P}^\epsilon (\cone) \big), 
\]
where $\mathbf{T}^\epsilon(\cone)$ is the $\epsilon$--Margulis tube about $\cone$ in $M$, the manifold $M_0$ is $M - \cone$ equipped with the complete hyperbolic structure whose conformal boundary agrees with that of $M$, and $\mathbf{P}^\epsilon (\cone)$ is the $\epsilon$--Margulis tube at the cusp corresponding to $\cone$. \qed
\end{theorem}

\noindent In fact, Proposition \ref{strainofcurvature}, Theorem \ref{meanvalue}, and Theorem \ref{L2bounddrilling} yield the following theorem.

\begin{theorem}[Brock--Bromberg, Corollary 6.10 of \cite{drilling}]\label{geodcurv} For any $\epsilon, \delta, C  > 0$ and $L > 1$, there exists an $\ell > 0$ so that if $\ell_{M_a}(\cone) < \ell$ then the following holds.
 Let $W$ be a subset of $M$, $\gamma(s)$ a smooth curve in $W$ and  $C(t)$ the geodesic curvature of $\gamma$ in the $g_t$--metric at $\gamma(0)$.  
If $W$ lies in the $\epsilon$--thick part $M_t^{\geq \epsilon}$ 
for all $t \in [t_0, a]$ and $C(0) \leq C$ then the identity map
\[
\mathrm{id} \co (W,g_a) \longrightarrow (W, g_{t_0}) 
\]
is $L$--bilipschitz and $| \, C(0) - C(a) \, | \leq \delta$. \qed
\end{theorem}

The map in Theorem \ref{drilling} is obtained by letting $W$ be the entire thick part and extending the map to the complement of a tubular neighborhood of the singular locus.  
This not only provides desirable bilipschitz maps, but will further allow us to control the change in the principal curvatures of certain surfaces under Dehn filling. 
We pause to formulate what we need.

\bigskip

\noindent Let $\F$ be an oriented Riemannian surface and let $M$ and $N$ be oriented Riemannian $3$--manifolds.  
Let 
\[
 \F \overset{f}{\to} M \overset{g}{\to} N 
 \]
be smooth with $f$ an immersion and $g$ an orientation preserving embedding. 
There are unique unit normal fields on $f(\F)$ and $g  f(\F)$ compatible with the orientations, and these are the unit normal fields that we use to define the normal curvatures of the two surfaces:
\[
\kappa_f(v) = \frac{\mathrm{II}_f(v,v)}{\| v \|_M^2} \quad \mathrm{and} \quad \kappa_{gf}(w) = \frac{\mathrm{II}_{gf}(w,w)}{\| w \|_N^2}
\]
where $\mathrm{II}_h$ is the second fundamental form of $h(\F)$.

We say that \textbf{the normal curvatures of $g  f (\F)$ are within $\delta$ of those of $f(\F)$} if for each $v$ in $T_p f(\F)$,
 \[
\big | \, \kappa_f(v) - \kappa_{gf}\big(Dg(v)\big) \big | < \delta.
 \]
Given a family of such maps
 \[
 \F_j \overset{f_j}{\to} M_j \overset{g_j}{\to} N_j \ ,
 \]
 we say that \textbf{the normal curvatures of $g_j  f_j(\F_j)$ tend to those of $f_j(\F_j)$} if they are within $\delta_j$ of each other with $\delta_j$ tending to zero. 
 We will need the following statement.

\begin{proposition}\label{curvdrilling}   For any $\epsilon, \delta, C > 0$ and $L > 1$, there exists an $\ell > 0$ so that if $\ell_{M_a}(\cone) < \ell$ the following holds.  
Let $\F$ be a smooth surface in $M = M_{t_0}$ such that $\F$ lies in the $\epsilon$--thick part $M_t^{\geq \epsilon}$  for all $t \in [t_0, a]$ and the principal curvatures of $\F$ are bounded by $C$. 
Then if $f_a \co M \to M_a$ is the $L$--bilipschitz map given by Theorem \ref{drilling},  the normal curvatures of $f_a(\F)$ are within $\delta$ of those of $\F$.

In particular, the principal curvatures of $f_a(\F)$ are within $\delta$ of those of $\F$. \qed
\end{proposition}

\subsection{Geometric inflexibility}
In addition to the geometric effect of cone deformations, 
we will need to control the effect of quasiconformal conjugacies of Kleinian groups deep in the hyperbolic manifolds. 
McMullen's Geometric Inflexibility Theorem says that in the convex core of a purely loxodromic Kleinian group, the effect of a quasiconformal deformation decays exponentially as one moves away from the boundary of the core. 
Precisely,


\begin{theorem}[McMullen, Theorem 2.11 of  \cite{renorm}]\label{inflexibility} Let $\Gamma$ be a finitely generated purely loxodromic Kleinian group and let $M_\Gamma = \HH^3/\Gamma$. 
A $K$--quasiconformal conjugacy of $\Gamma$ induces a map $f$ between hyperbolic manifolds whose bilipschitz constant at a point $p$ in the convex core $\mathcal{C}(M_\Gamma)$ satisfies
\[ \log \, \mathrm{bilip}(f)(p) \leq C e^{-A \dee(p, \, \partial \mathcal{C}(M_\Gamma))}\log K \] 
where $C$ is a universal constant and $A$ depends only on a lower bound on the injectivity radius at $p$ and an upper bound on the injectivity radius in $\mathcal{C}(M_\Gamma)$. \qed
\end{theorem}

Typically this theorem is used in the geometrically infinite setting, where one is free to move perpetually deeper into the core. 
We will be studying sequences of geometrically finite manifolds whose cores are getting deeper and deeper still while the dilatations of our deformations remain bounded. 
Unfortunately, our sequences will have convex cores of larger and still larger injectivity radii, and so McMullen's theorem will not be applicable. 
Fortunately, the following Geometric Inflexibility Theorem of Brock and Bromberg trades this bound for a bound on the area of the boundary and is sufficient for our purpose here.

\begin{theorem}[Brock--Bromberg \cite{bbinflex}]\label{inflex} Let $\Gamma$ be a finitely generated purely loxodromic Kleinian group and let $M_\Gamma = \HH^3/\Gamma$. 
A $\Gamma$--invariant Beltrami differential $\mu$ with $\| \mu \|_\smallinfinity < 1$ on $\widehat \C$ induces a map $f_\mu$ between hyperbolic manifolds whose bilipschitz constant at a point $p$ in the convex core $\mathcal{C}(M_\Gamma)$ satisfies
\[ 
\log \, \mathrm{bilip}(f_\mu)(p) \leq C e^{-A \dee(p, \partial \mathcal{C}(M_\Gamma))} \] 
where $C$ and $A$ depend only on $\| \mu \|_\smallinfinity$, a lower bound on the injectivity radius at $p$, and the area of $\partial \mathcal{C} (M_\Gamma)$. \qed
\end{theorem}

\noindent The map $f_\mu$ is constructed in a natural way that we now recall---see \cite{bbinflex} and Appendix B of \cite{renorm}.

By a visual averaging process, any vector field $v$ on $\widehat \C$ has a visual extension to a vector field $\mathrm{ex}(v)$ on $\HH^3$, and
any Beltrami differential $\mu$ on $\widehat \C$ has such an extension to a strain field $\mathrm{ex}(\mu)$.

\begin{theorem}[The Beltrami Isotopy] Let $f \co \widehat \C \to \widehat \C$ be a $K$--quasiconformal map normalized to fix $0$, $1$, and $\infinity$. 
Then there is a quasiconformal isotopy
\[
\varphi \co \widehat \C \times [0,1] \to \widehat \C
\]
fixing $0$, $1$, and $\infinity$ obtained by integrating a family $v_t$ of continuous $\frac{1}{2} \log K$--quasicon-formal vector fields, such that $\varphi_0 = \mathrm{id}$ and $\varphi_1 = f$.

The isotopy is natural in the sense that given $\gamma_0$ and $\gamma_1$ in $\PSL_2 \C$ such that 
$
\gamma_1 \circ f = f \circ \gamma_0 ,
$
there are interpolating M\"obius transformations $\gamma_t$ such that 
$
\gamma_t \circ \varphi_t = \varphi_t \circ \gamma_0 .
$ \qed
\end{theorem}

\noindent The Beltrami isotopy has a natural visual extension 
\[
\Phi \co \HH^3 \times [0,1] \to \HH^3
\]
that is the integral of the harmonic vector field $\mathrm{ex}(v_t)$. 
The maps $\Phi_t$ are bilipschitz, and letting $\eta_t$ denote the harmonic strain field $\mathrm{ex}(\overline{\partial} v_t)$, the pointwise bilipschitz constants satisfy
\[
\log \mathrm{bilip}(\Phi_t)(p) \leq \int_0^t \big\| \eta_s \big(\Phi_s(p) \big) \big\| \ \dee s \ .
\]

Taking our differential $\mu$ in Theorem \ref{inflex}, we begin with the normalized quasiconformal mapping with Beltrami coefficient $\mu$ and obtain an equivariant discussion, a $1$--parameter family of hyperbolic manifolds $M_t = (M, g_t)$, and bilipschitz $\Psi_t \co M_0 \to M_t$ obtained by pushing the $\Phi_t$ down to the quotient.
The harmonic strain fields $\mathrm{ex}(\overline{\partial} v_t)$ descend to those obtained from the $g_t$ by the construction in Section \ref{strain}, and at time one, this is the descendant of the strain field $\mathrm{ex}(\mu)$.

Our map $f_\mu$ equals $\Psi_1$, and its bilipschitz constant satisfies
\[
\log \mathrm{bilip}(f_\mu)(p) \leq \int_0^1 \big\| \eta_t \big(\Phi_t(p) \big) \big\| \ \dee \, t \ .
\]
The integral on the right is now estimated using

\begin{theorem}[Brock--Bromberg \cite{bbinflex}]\label{pointnorm} Let $\Gamma$ be a finitely generated purely loxodromic Kleinian group and $\mu$ a $\Gamma$--invariant Beltrami differential on $\widehat \C$. 
Let $K$ be defined by $\| \mu \|_\smallinfinity = (K-1)/(K+1)$ and let $b = \frac{1}{2}\log K$.
Let the harmonic strain field $\eta = \mathrm{ex}(\mu)$ be the visual extension of $\mu$.  
Then 
\[
\| \eta(p) \| + \| * D \eta(p) \| \leq 3A  b  \sqrt{\mathrm{area}(\partial \mathcal{C}(M))} \, e^{-\dee(p, \partial \mathcal{C}(M))} 
\]
where $p$ lies in the $\epsilon$--thick part of the convex core $\mathcal{C}(M)$ and $A$ depends only on $\epsilon$. \qed
\end{theorem}
\noindent Theorem \ref{quasiiso} and Proposition \ref{movinghulls} allow the transfer of this estimate to the desired estimate of $\log \mathrm{bilip}(f_\mu)(p)$, see \cite{bbinflex}.

As Theorem \ref{pointnorm} is not explicitly stated in \cite{bbinflex}, we derive it from the work there.
\begin{proof}[Derivation of Theorem \ref{pointnorm}]
The proof of Theorem 3.8 of \cite{bbinflex} and the fact that $* D \eta$ is a harmonic strain field establish the following theorem.
\begin{theorem}[Brock--Bromberg \cite{bbinflex}]\label{mod3.8}
Let $N$ be a complete hyperbolic $3$--manifold with compact boundary and let $\eta$ be a harmonic strain field on $N$. Then
\[
\| \eta(p) \| + \| * D \eta(p) \| \leq A  e^{-\dee(p, \partial N)} \sqrt{\int_N \| \eta \|^2 + \| * D \eta \|^2}
\]
where $p$ lies in the $\epsilon$--thick part of $N$ and $A$ depends only on $\epsilon$. 
\end{theorem}
The term 
\[
\int_N \| \eta \|^2 + \| * D \eta \|^2
\]
is bounded by Lemma 5.2 of \cite{bbinflex}:
\begin{lemma}[Brock--Bromberg, Lemma 5.2 of \cite{bbinflex}]\label{mod5.2} Let $N$ be a complete hyperbolic $3$--manifold such that $\pi_1(N)$ is finitely generated and assume that $N$ has no rank--one cusps. 
Let $\eta$ be a harmonic strain field on $N$ such that $\| \eta \|_\smallinfinity$ and $\| *D \eta \|_\smallinfinity$ are bounded by $B$. 
Then 
\[
\int_{\mathcal{C}(N)} \| \eta \|^2 + \| * D \eta \|^2 \leq \, B^2\, \mathrm{area}(\partial \mathcal{C}(M)).
\]
\end{lemma}
When $\eta$ is the strain field in Theorem \ref{pointnorm}, a theorem of H. M. Reimann \cite{reimann}, Part 1 of Theorem 5.1 of \cite{bbinflex}, says that $\| \eta \|_\smallinfinity$ and $\| *D \eta \|_\smallinfinity$ are bounded by $3b$.
Theorem \ref{pointnorm} now follows by letting $N = \mathcal{C}(M)$ in Theorem \ref{mod3.8} and Lemma \ref{mod5.2}.
\end{proof}

The estimate in Theorem \ref{pointnorm} provides the following Proposition, thanks to Proposition \ref{strainofcurvature}.

\begin{proposition}\label{curvinflex}  For any $\epsilon, \delta, C > 0$ and $0 < k < 1$, there is a $d$ such that the following holds. 
Let $\mu$ be a Beltrami differential for a finitely generated purely loxodromic Kleinian group $\Gamma$ with $\| \mu \|_\smallinfinity = k$ and let $M = M_\Gamma$. 
Let $\F$ be a smooth surface in $M^{\geq \epsilon} \cap \mathcal{C}(M)$  at a distance at least $d$ from $\partial \mathcal{C}(M)$, and suppose that the principal curvatures of $\F$ are bounded by $C$. 
Then if $f_\mu \co M \to M_\mu$ is the map in Theorem \ref{inflex},  the normal curvatures of $f_\mu(\F)$ are within $\delta$ of those of $\F$. 

In particular, the principal curvatures of $f_\mu(\F)$ are within $\delta$ of those of $\F$. \qed 
\end{proposition}

\subsection{Useful surfaces}\label{usefulsurfaces}

Let $\F$ be a smooth surface in $\HH^3$.  
Pick a point $p$ in $\F$ and normalize so that $p$ sits at $(0,0,1)$ in the upper half-space model with unit normal $-\frac{\partial}{\partial z}$ 
such that the principal directions at $p$ are $(1,0,0)$ and $(0,1,0)$. 
Let $\Pi$ be the normal projection of $\F$ to $\C$. 
Picking an orthonormal basis for $T_p\F$ along its principal directions and the usual basis for $T_0 \C$, the derivative of $\Pi$ at $p$ is given by the matrix
\begin{equation}\label{derivative}
D \Pi_p = \begin{pmatrix} \frac{1 + \kappa_1}{2}  &   0 \\
  0 &  \frac{1 + \kappa_2}{2}
\end{pmatrix}
\end{equation}
where the $\kappa_i$ are the principal curvatures of $\F$---the reader will find it a nice exercise in plane geometry to verify the formula using the fact that, with our normalization, the hyperbolic curvature at $(0,0,1)$ in a vertical hyperbolic plane is one plus the Euclidean curvature (a proof of \textit{this} fact may be found in \cite{anderson}).

\subsubsection{The Epstein surface}

Given a quasifuchsian manifold $\QF(X,Y)$ with convex core $\mathcal{C}(X,Y)$ possessing boundaries $\mathcal{O}_X$ and $\mathcal{I}_Y$ facing $X$ and $Y$ respectively, we let $\mathcal{E}_X$ be the component of $\QF(X,Y) - \mathcal{C}(X,Y)$ facing $X$, $\mathcal{E}_Y$ the component facing $Y$.

\begin{lemma}\label{oursurface} Let $\mathcal{K}$ be a compact set in $\T(\overline{\Sigma})$. 
There is a number $d > 0$ such that for any quasifuchsian manifold $\QF(X,Y)$ in $\QF(\T(\Sigma), \mathcal{K})$ there is a smooth convex surface in $\mathcal{E}_Y$ that lies in the $d$--neighborhood of $\mathcal{I}_Y$ and whose principal curvatures are greater than $\frac{1}{2}$ and no more than $2$. 
\end{lemma}

Given a conformal metric $m$ on a simply connected domain $\Omega$ in $\widehat \C$, there is a unique map $E \co (\Omega, m) \to \HH^3$ with the property that the visual metric at $E(z)$ agrees with $m$ at $z$ in $\Omega$. 
This map was discovered by C. Epstein \cite{epstein}, and the following information concerning it may be found in the work of Epstein \cite{epstein} and C. G. Anderson \cite{anderson}.
A nice summary of the information we need here may be found in \cite{brombergJAMS}.

We will only be interested in scalar multiples of the Poincar\'e metric $\rho$ on $\Omega$, and we write $E_t \co (\Omega, e^t \rho) \to \HH^3$.  

We define the norm
\[
\| \, \cdot \, \|_{\smallinfinity,\rho} = \big \| \, \cdot \, (1- |z|^2)^{-2} \big\|_{L^\scriptinfinity (\ELL)} 
\]

Let $f_\Omega$ be a Riemann mapping carrying the unit disk $\ELL$ to $\Omega$.  
By Theorem 7.5 of \cite{epstein}, there is a $T$ depending only on an upper bound for the weighted $L^\smallinfinity$--norm $\| \, Sf_\Omega \, \|_{\smallinfinity,\rho}$ of the Schwarzian derivative $Sf_\Omega$ of $f_\Omega$ such that $E_s$ is an embedding with strictly convex image whenever $s \geq T$---the surface is curving \textit{away} from $\Omega$. 
Thanks to a theorem of Z. Nehari and W. Kraus, the norm $\| Sf_\Omega \|_{\smallinfinity,\rho}$ is bounded by six, see Section 5.4 of \cite{gardiner}.  
So we may take $T$ universal. 

We now restrict our attention to $\Omega$ that arise as a component of the domain of discontinuity of a quasifuchsian group $\QF(X,Y)$---we will always think of $\Omega$ as the domain uniformizing $Y$.  The principal curvatures of the surface $E_t(\Omega)$ tend to one independently of $\Omega$ as $t$ tends to infinity, see Section 3 of \cite{anderson} or Proposition 6.3 of \cite{brombergJAMS}, and so we assume as we may that for all $s$ greater than or equal to $T$, these curvatures are at least $1/2$ and no more than $2$.

All of the above works $\Gamma_{\QF(X,Y)}$--equivariantly 
and so, for any $\QF(X,Y)$, we have a smoothly embedded convex copy $E_T(Y)$ of $Y$ in $\mathcal{E}_Y$ with principal curvatures trapped between $1/2$ and $2$. 
The map $E_T \co Y \to \QF(X,Y)$ is Lipschitz with constant 
\[
\| S f_\Omega \|_{\smallinfinity, \rho} + \frac{1}{2}(e^T + e^{-T})  \leq 6 + \frac{1}{2}(e^T + e^{-T}),
\]
 see \cite{anderson}.

We are now ready for the proof of the lemma.
\begin{proof}[Proof of Lemma \ref{oursurface}] All that remains to be shown is that the Hausdorff distance between $\mathcal{I}_Y$ and $E_T(Y)$ depends only on $T$ and the compact set $\mathcal{K}$.  
Since the hyperbolic surfaces in $\mathcal{K}$ have uniformly bounded diameter, and the Lipschitz constant of $E_T$ is less than $6 + \frac{1}{2}(e^T + e^{-T})$, we obtain a bound $A$ on the diameter of the surfaces $E_T(Y)$ that only depends on $T$ and $\mathcal{K}$. 

Now, there is a universal lower bound $D$ to the diameter of any complete hyperbolic surface, and so the diameter of the boundary of the $R$--neighborhood $N_R(\mathcal{I}_Y)$ of $\mathcal{I}_Y$ in $\mathcal{E}_Y$ is at least $D\cosh R$. 
If $E_T(Y)$ is entirely contained in the complement of $N_R(\mathcal{I}_Y)$, its diameter is at least $D\cosh R$, as the projection $E_T(Y) \to \partial N_R(\mathcal{I}_Y)$  is $1$--Lipschitz and has degree one.  
So $E_T(Y)$ must intersect the $\mathrm{arccosh}(2A/D)$--neighborhood of $\mathcal{I}_Y$.
Letting $d = A + \mathrm{arccosh}(2A/D)$ completes the proof.
\end{proof}

\section{Focusing}\label{focusing} 

Given a totally geodesic surface $\Sigma$ in a hyperbolic $3$--manifold $M$, a \textbf{metric collar} of \textbf{depth} $d$ about $\Sigma$ is a $d$--neighborhood of $\Sigma$ in $M$ homeomorphic to $\Sigma \times I$. 

The following theorem says that the skinning map of a manifold with a deep collar about its boundary is highly contracting on a large compact set.
\begin{theorem}[Focused Manifolds]\label{collar} Let $S$ be a closed hyperbolic surface. 
Let $\epsilon , \delta, R > 0$. 
There is a constant $d > 0$ such that if $M$ is a hyperbolic $3$--manifold with totally geodesic boundary $\Sigma$ homeomorphic to $S$ with $\mathrm{injrad}(\Sigma) > \delta$ possessing a collar of depth $d$, then $\diam \big( \sigma_M\big( B(\Sigma,R) \big) \big) < \epsilon$.
\end{theorem}
\begin{proof} Fix $\delta, R > 0$ and let $\{M_d\}_{d = 1}^\smallinfinity$ be \textit{any} sequence of hyperbolic manifolds with totally geodesic boundary such that  $\Sigma_d = \partial M_d$ is homeomorphic to $S$, has injectivity radius at least $\delta$, and has a metric collar of depth $d$ in $M_d$. 
We write $\sigma_d = \sigma_{M_d}$.
We show that the diameter of $\sigma_d\big( B(\Sigma_d,R) \big)$ is less than $\epsilon$ for all sufficiently large $d$.

We consider $\pi_1(\dot M_d)$ with basepoint in $\Sigma_d$ as a Kleinian group $\Gamma_d$ normalized so that the complement $\U$ of the closed unit disk $\overline{\ELL}$ in $\widehat \C$ is the component of the domain of discontinuity $\Omega_{\Gamma_d}$ uniformizing $\Sigma_d$ with uniformizing Fuchsian group $\Gamma_{\Sigma_d}$.

Let  $\mu = \mu_d$ be a Beltrami differential of norm $k \leq (e^{2R} - 1)/(e^{2R} + 1)$ for $\Gamma_{\Sigma_d}$ in $\U$, so that $[\mu]$ lies in $B(\Sigma_d, R)$.

Let $\F_d$ be the component of the boundary of the $d/2$--neighborhood of $\Sigma_d$ contained in the convex core of $M_d$.  
Since $\Sigma_d$ is totally geodesic, this is a smooth surface.    
Moreover, it is easy to see that the principal curvatures of the $\F_d$  tend to one as $d$ grows.  
Since the injectivity radius of $\Sigma_d$ is at least $\delta$, we may assume that the injectivity radius of $M_d$ is at least $\delta$ at $\F_d$, by choosing $d$ larger than $\delta$.

The differential $\mu$ induces a bilipschitz map $f_{d,\mu} \co M_d \to M_{d,\mu}$, and we write $\F_{d,\mu} = f_{d,\mu}(\F)$.
As $d$ grows, the bilipschitz constants of the $f_{d,\mu}$ near $\F_d$ tend to one in a manner depending only on $d$, $R$, $\delta$, and the topological type of $S$, by Theorem \ref{inflex}.

Since the norm of $\mu$ is bounded, Proposition \ref{curvinflex} tells us that for any small $\eta > 0$, there is a large $d$ such that the principal curvatures of the $\F_{d,\mu}$ are within $\eta$ of those of $\F_d$. 
As the principal curvatures of $\F_d$ are universally bounded away from zero, the $\F_{d,\mu}$ are strictly convex for large $d$, and
so there are diffeomorphic normal projections 
\[
\Pi_{d,\mu} \co \sigma_d([\mu]) \to \F_{d,\mu}
\quad
\mathrm{and}
\quad
\Pi_d \co  \Sigma_d \to \F_d .
\]

Consider the composition
\[
F_{d,\mu} = \Pi_{d,\mu}^{-1} \circ f_{d,\mu} \circ \Pi_d 
\co \Sigma_d
\longrightarrow \sigma_d([\mu]) 
\]
and its derivative
\[
DF_{d,\mu} = D\Pi_{d,\mu}^{-1} \circ Df_{d,\mu} \circ D\Pi_d 
\]
at a point $p$ in $\Sigma_d$ (we suppress $p$ in the notation).

By the discussion in Section \ref{usefulsurfaces}, we may choose coordinates so that
\[
DF_{d,\mu} =
A^{-1}
\begin{pmatrix} \frac{1 + \kappa_1(\F_{d,\mu})}{2}  &   0 \\
  0 &  \frac{1 + \kappa_2(\F_{d,\mu})}{2}
\end{pmatrix}
A
\begin{pmatrix} w'  &   x' \\
  y' & z' 
\end{pmatrix}
\begin{pmatrix} \frac{2}{1 + \kappa_1(\F_d)}  &   0 \\
  0 &  \frac{2}{1 + \kappa_2(\F_d)}
\end{pmatrix},
\]
where $A$ is the rotation carrying the principal directions of $\F_{d,\mu}$ to $(1,0)$ and $(0,1)$.


To estimate the Teichm\"uller distance between $\Sigma_d$ and $\sigma_d([\mu])$, we would like to understand the ``norm" of the image of a vector under $DF_{d,\mu}$.
Unfortunately, our skinning surfaces do not come to us equipped with tractable Riemannian metrics, and so we must be careful when using information about $DF_{d,\mu}$ to estimate the Teichm\"uller distance between them.
This is achieved as follows.

Note that with our normalization the hyperbolic inner products on the tangent spaces to our surfaces $\F_d$ and $\F_{d,\mu}$ agree with the standard Euclidean ones, as we have chosen charts so that they are tangent to the horosphere of height one.
Our coordinate patches for $\Sigma_d$ and $\sigma_d([\mu])$ lie centered at $0$ in $\C$, and we equip their tangent spaces at zero with the standard Euclidean inner products, which are compatible with the conformal structures on $\Sigma_d$ and $\sigma_d([\mu])$.

Now, if $\lambda_+ \geq \lambda_-$ are the maximum and minimum values of $\| \, DF_{d,\mu} v \, \|$ as $v$ ranges over the unit circle in $T_p \Sigma_d$, then the dilatation of $F_{d,\mu}$ at $p$ is no more than $\lambda_+/\lambda_-$, see the first chapter of \cite{ahlfors}.
It follows that if the linear map $DF_{d,\mu}$ is $\ell$--bilipschitz, the map $F_{d,\mu}$ is $\ell^2$--quasiconformal at $p$.


Writing
\[
\begin{pmatrix} w  &   x \\
  y & z 
\end{pmatrix}
=
A
\begin{pmatrix} w'  &  x' \\
  y' & z' 
\end{pmatrix}
\]
and multiplying, we have
\[
DF_{d,\mu} =
A^{-1}B =
A^{-1}
\begin{pmatrix} w\,  \frac{1 + \kappa_1(\F_{d,\mu})}{1 + \kappa_1(\F_d)}  
&   x \,  \frac{1 + \kappa_1(\F_{d,\mu})}{1 + \kappa_2(\F_d)} 
\\[7pt]
  y \,  \frac{1 + \kappa_2(\F_{d,\mu})}{1 + \kappa_1(\F_d)} 
  & z \, \frac{1 + \kappa_2(\F_{d,\mu})}{1 + \kappa_2(\F_d)} 
\end{pmatrix}.
\]
The numbers $w$, $x$, $y$, and $z$ are bounded above in absolute value by a constant depending only on the bilipschitz constant of $f_{d,\mu}$, which is bounded independent of $\mu$.

Let $\eta > 0$ be very small.  
For large $d$, the principal curvatures of the surfaces involved are all within $\eta$ of one.
So the bilipschitz constant of the linear map $B$, and hence that of $A^{-1}B$, is very close to that of $Df_{d,\mu}$ in a manner depending only on $d$, $R$, $\delta$, and the topological type of $S$.
Since the $f_{d,\mu}$ are bilipschitz with constants tending to one as $d$ grows, the same is true for the linear maps $DF_{d,\mu}$.
So the $F_{d,\mu}$ are quasiconformal with constants tending to one as $d$ grows, and $\Sigma_d$ and $\sigma_d([\mu])$ are close in the Teichm\"uller metric independent of $\mu$.
\end{proof}

\subsection{Ubiquity of focus}\label{ubiquity}

To prove our Quasifuchsian Density Theorem in Section \ref{quasifuchsiandensitysection}, and the Capping Theorem in Section \ref{capsection}, we will need the following theorem.

\begin{theorem}\label{small} Let $X$ be a closed connected hyperbolic surface, let $B(X,R)$ be the ball of radius $R$ about $X$ in $\T(X)$, and let $\epsilon > 0$.  
Then there is a convex cocompact hyperbolic manifold $M$ with totally geodesic boundary $\Sigma$ homeomorphic to $X$ satisfying
\[
\dee(\Sigma, X) < \epsilon
\quad
and
\quad
\diam \big( \sigma_M\big( B(X,R) \big) \big) < \epsilon.
\]
\end{theorem}

Let $\mathfrak{C}_d$ be the class of convex cocompact hyperbolic manifolds with totally geodesic boundary whose boundaries have collars of depth $d$.
By Theorem \ref{collar}, to prove the theorem we need only prove the following mild generalization of Fujii and Soma's theorem that the set of totally geodesic boundaries of hyperbolic manifolds is dense in the Teichm\"uller space \cite{fujiisoma}.
\begin{theorem}\label{ddensity} For any $d \geq 0$, the set of hyperbolic surfaces in $\T(S)$ appearing as totally geodesic boundaries of hyperbolic $3$--manifolds in $\mathfrak{C}_d$ whose boundaries are homeomorphic to $S$ is dense. 
\end{theorem}

A circle packing on a hyperbolic surface is a collection of circles with disjoint interiors whose union has the property that all of the complementary regions are curvilinear triangles.  
A hyperbolic structure on a surface is a \textbf{circle packing hyperbolic structure} if it admits a circle packing. 

We have the following theorem of R. Brooks. 

\begin{theorem}[Brooks \cite{brooks,brooksschottky}] The set of circle packing hyperbolic structures on a surface is dense in the Teichm\"uller space. \qed
\end{theorem}

\noindent Define the \textbf{radius} of a circle packing to be the maximum of the radii of its circles.
The work in \cite{brooksschottky} and \cite{brooks} in fact proves the following theorem.

\begin{theorem}[Brooks]\label{strongbrooks} Let $\epsilon > 0$.  
Then the set of hyperbolic structures on a surface admitting circle packings of radius at most $\epsilon$ is dense in the Teichm\"uller space.
\end{theorem}
\begin{proof}[Sketch (see the proof of Theorem 3 of \cite{brooks}).] Let $X$ be a hyperbolic surface.
Any configuration of circles in $X$ may be extended to a configuration whose complementary regions are curvilinear triangles and rectangles.  
We choose such a configuration consisting of circles of radius $\epsilon/2$ whose complementary regions have diameter at most $\epsilon/2$.

Consider the group $\Gamma$ of M\"obius transformations generated by the Fuchsian group uniformizing $X$ and all of the reflections in the circles of our configuration. 
As in the proof of Theorem 3 of \cite{brooks}, we may choose an arbitrarily small quasiconformal deformation of $\Gamma$ carrying our configuration of circles to a configuration that fits into a circle packing.
Since the deformation may be taken as small as desired,  we may take the radius of the resulting circle packing less than $\epsilon$. 
The fact that this deformation induces a small quasiconformal deformation of $X$ completes the proof.
\end{proof}

\begin{proof}[Proof of Theorem \ref{ddensity}] We begin the proof as Fujii and Soma do, but use the Bounded Image Theorem (which Fujii and Soma avoid) and McMullen's Geometric Inflexibility Theorem to complete it.

Fix a manifold $B$ with totally geodesic boundary homeomorphic to $S$ for later use,
and let $\mathcal{K}$ be a compact set in $\T(\overline{S})$ containing the image of $\sigma_{B}$.

Let $d, \delta > 0$.

We will show that for any surface $X$ in $\T(S)$ there is a manifold in $\mathfrak{C}_d$ whose totally geodesic boundary is within $\delta$ of  $X$.

Let $\Sigma$ be a hyperbolic surface equipped with a circle packing $\mathcal{C}$ of radius $\epsilon$.  By Theorem \ref{strongbrooks}, these are dense in $\T(S)$ for any $\epsilon$. 
We prove that if $\epsilon$ is small enough, then there is a convex cocompact manifold with geodesic boundary within $\delta/2$ of $\Sigma$ that possesses a collar of depth $d$.
It follows that \textit{any} point in $\T(S)$ is within $\delta$ of the boundary of such a manifold.

By inscribing each complementary triangle in a circle, we obtain a dual configuration of circles $\mathcal{D}$. 
Consider the group $\Delta$ of M\"obius transformations generated by $\Gamma_\Sigma$ and the reflections  in the circles in $\mathcal{C} \cup \mathcal{D}$. 
This group has a finite index torsion free subgroup whose Kleinian manifold has four boundary components, each conformally or anticonformally equivalent to $\Sigma$.  
Moreover, the corresponding convex core has totally geodesic boundary.

This manifold may be constructed directly by ``scalloping" the circle packing, as we now discuss.  
Lift the configurations $\mathcal{C}$ and $\mathcal{D}$ to the unit disk in $\widehat \C = \partial_\smallinfinity \HH^3$. 
Each circle is the ideal boundary of a totally geodesic hyperplane in $\HH^3$, and for each such circle we excise from $\HH^3$ the open halfspace incident to the interior of the circle.  
We also excise the halfspace corresponding to the complement of the unit disk.  
What remains is a convex subset of $\HH^3$ whose boundary is a union of totally geodesic polygons meeting at right angles, together with a single copy of $\HH^2$.

Let $\mathcal{O}$ be the quotient of this convex set by $\Gamma_\Sigma$.  
The set $\partial \mathcal{O} - \Sigma$ may be decomposed into two families, a family of triangles $T$ arising from the triangles complementary to $\mathcal{C}$, and the remaining polygons. 
Now, we double $\mathcal{O}$ along $T$, and then double this manifold along the double of the remaining polygons.  
The resulting manifold has four totally geodesic boundary components as described, and we glue two of them together to obtain a manifold $M$ with $\partial M = \Sigma \cup \overline{\Sigma}$.

Notice that since the radius of $\mathcal{C}$ is less than $\epsilon$, the boundary of $M$ has a very large collar, as all of the totally geodesic planes in the universal cover must be far away from the one whose boundary is the unit circle.

The manifold $M$ has cusps, and these may be filled using the Hyperbolic Dehn Filling Theorem so that the totally geodesic boundary of the resulting manifold is very close to that of $M$ and so that the filled manifold still possesses a large collar about its boundary.
More precisely, consider the double $DM$ of $M$ along its boundary.
If we perform higher and higher hyperbolic Dehn fillings equivariantly with respect to the natural involution of $DM$, the resulting manifolds converge geometrically to $DM$, see Section E.1 of \cite{bp}.  
Moreover, these fillings are themselves doubles, and so the copies of $\partial M$ in these manifolds are isotopic to totally geodesic surfaces.  
The geometric convergence implies that these Fuchsian groups are converging algebraically to the corresponding Fuchsian group in $DM$, and so the totally geodesic surfaces in these fillings are eventually within $\delta/4$ of $\Sigma$ in $\T(S)$.
The geometric convergence also guarantees a large collar about these surfaces, and cutting open again provides the desired manifolds.

We continue the argument with one of these filled manifolds $M'$ with boundary $\Sigma' \cup \overline{\Sigma}'$.

For definiteness, we say that the depth of the collar about the boundary of $M'$ is $c$.
By the above, $c$ tends to infinity as $\epsilon$ goes to zero.

Now attach $B$ to the copy of $\overline{\Sigma}'$ in $\partial M'$ to obtain a manifold $N$.
The manifold $N$ is irreducible, acylindrical, and atoroidal, as is seen by considering potential essential spheres, annuli, and tori in $N$, which, if there, would result in a sphere, annulus, or torus in $M'$ or $B$---see Lemma 2.1 of \cite{myersexcellent}. 
By Thurston's Hyperbolization Theorem for Haken manifolds, the manifold $N$ admits a hyperbolic structure with totally geodesic boundary, and we equip $N$ with this structure. 

We claim that for large $c$, the totally geodesic boundary of $N$ is within $\delta/4$ of  $\Sigma'$ (and hence within $\delta$ of $\Sigma$) with a collar of depth $d$.

To see this, consider the double $DN$ of $N$.  This manifold is obtained from a double $DM'$ of $M'$ by capping off with $B \sqcup \overline{B}$.  
We have a map
\[
\iota \co \GF(DN) \to \GF(DM') \cong \T(S) \times \T(\overline{S})
\]
where, of course, the domain is a single point, by G. D. Mostow's Rigidity Theorem \cite{mostow}. 
The image of $\iota$ lies in $\overline{\mathcal{K}} \times \mathcal{K}$.

In $DM'$, the surface $\Sigma_0$ that will be the boundary of $N$ is at a distance $c$ from $\partial DM'$. 
By McMullen's Geometric Inflexibility Theorem, the bilipschitz constant of the map between the hyperbolic structure on $DM'$ with totally geodesic boundary and the image of $\iota$ is very close to one on a large neighborhood of $\Sigma_0$ when $c$ is very large, as the diameter of $\mathcal{K}$ is fixed. 
We conclude that the totally geodesic boundary of $N$ has a collar of depth $d$ provided $c$ is large enough.
We also see that as $c$ tends to infinity, the Fuchsian group in $\Gamma_{DN}$ corresponding to $\Sigma_0$ is tending to the Fuchsian group $\Gamma_{\Sigma_0} = \Gamma_{\Sigma'}$ in the algebraic topology.  
We conclude that the totally geodesic boundary of $N$ is within $\delta/4$ of $\Sigma'$, and hence within $\delta/2$ of $\Sigma$.
\end{proof}

\subsection{Quasifuchsian density}\label{quasifuchsiandensitysection}

A theorem of Brooks' \cite{brooks} says that any convex cocompact Kleinian group embeds in a cocompact Kleinian group after an arbitrarily small quasiconformal deformation.  
In case the Kleinian group is quasifuchsian, Brooks' argument produces a cocompact Kleinian group in which the quasifuchsian surface is nonseparating. 
We have the following.
\begin{theorem}[Quasifuchsian Density] The subset of  $\QF(S)$ consisting of quasifuchsian groups that appear as the only boundary subgroup of a geometrically finite hyperbolic $3$--manifold is dense. 
\end{theorem}
\begin{proof} Let $\QF(X,\overline{Y})$ be a quasifuchsian group in $\QF(S)$ and let $\epsilon>0$.  
Let $R$ be such that $X$ lies in $B(Y,R)$. 
By Theorem \ref{small}, there is a manifold $M$ with totally geodesic boundary $Z$ within $\epsilon/2$ of $Y$ such that $\sigma_M \big(B(Y,R) \big)$ has diameter less than $\epsilon/2$. 
So $\sigma(X)$ is within $\epsilon$ of $\overline{Y}$. 
\end{proof}

\subsection{Capping off}\label{capsection}

The following theorem is useful in promoting properties of skinning maps of manifolds with many boundary components to those of manifolds with connected boundary.

\begin{theorem}[Capping Theorem]\label{cap} Let $M$ be a hyperbolic $3$--manifold with totally geodesic boundary $\Sigma_0 \sqcup \Sigma_1 \sqcup\ldots \sqcup \Sigma_m$.  
Let $\epsilon>0$. 
Then $(M, \Sigma_0)$ embeds as a pair into a hyperbolic manifold $(M_\epsilon, \partial M_\epsilon)$ with totally geodesic boundary homeomorphic to $\Sigma_0$ with the property that 
\[
\dee_{\T(\overline{\Sigma}_0)} \big(   \sigma_{M,\Sigma_0}(X),    \sigma_{M_\epsilon}(X)   \big) < \epsilon
\]
for all $X$ in $\T(\Sigma_0)$, where $\sigma_{M,\Sigma_0}$ is the skinning map of $M$ relative to $\Sigma_1 \sqcup\ldots \sqcup \Sigma_m$.
\end{theorem}
\begin{proof} By Theorem \ref{small}, for any $\epsilon$ and $R$ greater than zero there are hyperbolic $3$--manifolds $N_1, \ldots, N_{m}$ with totally geodesic boundaries $Y_1, \ldots, Y_{m}$ whose skinning maps satisfy
\begin{equation}\label{diamestimate}
\mathrm{diam}\Big(\sigma_{N_j}\big(B(\Sigma_j,R)\big)\Big) < \frac{1}{4}\log(1+\epsilon)^{2/3}
\end{equation}
and
\begin{equation}\label{boundaryestimate}
\dee\Big(\overline{\Sigma}_j, \sigma_{N_j}\big(B(\Sigma_j,R)\big)\Big) < \frac{1}{4}\log(1+\epsilon)^{2/3}.
\end{equation}

Now, let $R = \diam\big(\sigma_{M} \big(\T(S) \big)\big)$ and form a $3$--manifold $M_\epsilon$ by attaching $N_j$ as above to the $\Sigma_j$ in $\partial \dot M$.
We have the maps 
\[
\iota^* \co \GF(M_\epsilon) 
\longrightarrow \GF(\dot M) \cong \T(\Sigma_0) \times \T(\Sigma_1) \times \cdots \times \T(\Sigma_{m})
\]
and
\[
Q_0^* \co \GF(\dot M) \to \QF(\Sigma_0)
\]
induced by inclusion.
By \eqref{diamestimate}, \eqref{boundaryestimate}, Theorem \ref{quasiiso}, and the triangle inequality, for any $X$, the hyperbolic metrics 
$\iota^*(X)$ and $(X,\Sigma_1,\ldots,\Sigma_{m})$ are $(1+\epsilon)$--bilipschitz. 
In particular, the quasifuchsian manifolds 
\[
Q_0^*\iota^*(X)= \QF \big(X, \sigma_{M_\epsilon}(X) \big)
\quad
\mathrm{and}
\quad
\QF \big(X, \sigma_{M,\Sigma_0}(X) \big)
\]
are $(1+ \epsilon)$--bilipschitz. 
Theorem \ref{quasiiso} completes the proof.
\end{proof}

\section{Filling}\label{fillingsection}

The following theorem allows us to perform hyperbolic Dehn filling without terribly disturbing the skinning map.

\begin{theorem}[Filling Theorem]\label{fill} Let $M$ be a finite volume hyperbolic manifold with nonempty closed totally geodesic boundary and let $\mathfrak{e} > 0$. 
There is an $\hbar = \hbar_\mathfrak{e}(M) > 0$ such that if the normalized length of each component of $\alpha$ in $M$ is at least $\hbar$, then
\[
\dee_{\T(\overline{\partial M})}\big(\sigma_M(X), \sigma_{M(\alpha)}(X) \big) < \mathfrak{e}
\]
for all $X$ in $\T(\partial M)$.\end{theorem}
\begin{proof} For simplicity, we assume that $\partial M$ is connected.
We let $\ell(\alpha)$ be the minimum of the normalized lengths of the components of $\alpha$.

Let $M_X$ be $M^\circ$ equipped with the hyperbolic metric corresponding to $X$ in $\T(\partial M)$.
Let $\mathbf{P}_X^\epsilon$ denote the $\epsilon$--Margulis tube at the cusps of $M_X$, 
and let $\mathbf{T}_X^\epsilon(\alpha)$ denote the union of the $\epsilon$--Margulis tubes at the cusps of $M_X(\alpha)$ and about the components of the geodesic core $\cone$ of the union of filling tori in $M_X(\alpha)$.

Let $\mathcal{I}_X$ be the component of the boundary of the convex core of $\QF \big(X,\sigma_M(X) \big)$ facing $\sigma_M(X)$, and let $\F_X \subset \QF \big(X, \sigma_M(X) \big)$ be the smooth surface facing $\sigma_M(X)$ in the $d$--neighborhood of $\mathcal{I}_X$ given by the Bounded Image Theorem and Lemma \ref{oursurface}. 
We have the restriction 
\[
f_X \co \F_X \to M_X
\]
of the covering map $\QF \big(X, \sigma_M(X) \big) \to M_X$, and we claim that there is an $\epsilon$ such that, for all $X$, the surface $f_X(\F_X)$ lies in the $\epsilon$--thick part of $M_X$. 

By the Bounded Image Theorem, there is an  $\epsilon_0  >0 $ such that the image of $\sigma_M$ is $\epsilon_0$--thick.
By a theorem of Sullivan \cite{sullivan}, there is a universal $K$ for which $\mathcal{I}_X$ and $\sigma_M(X)$ are $K$--bilipschitz, see \cite{epsteinmarden}.
So the injectivity radius of $\mathcal{I}_X$ is at least $\epsilon_0/K$.
Taking $\epsilon_1$ less than the constant $\delta_{\epsilon_0/2K} = \delta_{\epsilon_0/2K}(\partial M)$ given by Lemma \ref{brooksmatelskilemma} will ensure that $\mathcal{I}_X$ lies in the $\epsilon_1$--thick part of $M_X$.
(We also assume that $\epsilon_1$ is less than the $\delta_{\epsilon_3}(\partial M)$ given by Lemma \ref{brooksmatelskilemma}, a fact we use in the claim below.)
Since $\F_X$ lies in the $d$--neighborhood of $\mathcal{I}_X$, this provides the desired $\epsilon < \epsilon_1 < \delta_{\epsilon_0/2K}$, which we now fix.


Let $h$ be the constant given by Theorems \ref{universal} and \ref{flatbounded} such that the components of the geodesic core $\cone$ have length less than $\epsilon$ when $\ell(\alpha) \geq h$, and assume that $\ell(\alpha)$ is this large.


\bigskip \noindent
Now let $L >1$ and $\eta > 0$.
Theorem \ref{drilling} and Proposition \ref{curvdrilling} provide an $H \geq h$ and $L$--bilipschitz maps 
\[
F_{X,\alpha} \co M_X - \mathbf{P}_X^\epsilon 
\longrightarrow M_X(\alpha) - \mathbf{T}_X^\epsilon(\alpha)
\]
when $\ell(\alpha) \geq H$ such that the principal curvatures of the $\mathfrak{G}_{X,\alpha} = F_{X,\alpha}\big(f_X(\F_X)\big)$ are within $\eta$ of those of $f_X(\F_X)$. 

The map 
\[
F_{X,\alpha} \circ f_X \co \F_X \longrightarrow M_X(\alpha) - \mathbf{T}_X^\epsilon(\alpha)
\]
lifts to a map
\[
G_{X,\alpha} \co \F_X \longrightarrow \QF \big(X, \sigma_{M(\alpha)}(X) \big),
\]
whose image we call
\[
\F_{X,\alpha} = G_{X,\alpha}(\F_X).
\]


\begin{claim} There is an $\hbar \geq H$ depending only on $H$ and the diameter and thickness of $\sigma_M$ such that $G_{X,\alpha}$ is an embedding provided $\ell(\alpha) \geq \hbar$.
\end{claim}
\begin{proof}[Proof of claim]
If $\gamma$ is an element of $\pi_1 \big(M(\alpha) \big)$ corresponding to a component of $\cone$, then
no power of $\gamma$ is conjugate into $\pi_1\big(\partial M(\alpha) \big)$.
To see this, consider the double $DM$.
Filling equivariantly with respect to the involution, we find that $M(\alpha)$ admits a hyperbolic structure $\mathcal{C} \big(M_Y(\alpha) \big)$ with totally geodesic boundary $Y$, and that the geodesic cores of the filling tori lie in the interior of $\mathcal{C} \big(M_Y(\alpha) \big)$.
Since $\partial \mathcal{C} \big(M_Y(\alpha) \big)$ is totally geodesic, no power of $\gamma$ is conjugate into $\pi_1\big(\partial M(\alpha) \big)$.

Let 
\[
b \co \QF\big(X, \sigma_{M}(X)\big) \to M_X
\]
be the covering corresponding to $\pi_1(\partial \dot M_X)$ and let
\[
b_\alpha \co \QF\big(X, \sigma_{M(\alpha)}(X)\big) \to M_X(\alpha)
\]
be the covering corresponding to $\pi_1\big(\partial \dot M_X(\alpha)\big)$.
For any $\omega > 0$, write
\[
\mathbf{B}^\omega = b^{-1}\big(\mathbf{P}_X^\omega\big)
\quad
\mathrm{and}
\quad
\mathbf{B}_\alpha^\omega = b_\alpha^{-1}\big(\mathbf{T}_X^\omega(\alpha)\big).
\]

Let $\mathcal{C}_\alpha$ be the convex core of $\QF\big(X, \sigma_{M(\alpha)}(X)\big)$,
and let $\mathcal{I}_\alpha$ be the component of $\partial \mathcal{C}_\alpha$ facing $\sigma_{M(\alpha)}(X)$.
By our choice of $\epsilon$ (smaller than the $\delta_{\epsilon_3}$ in Lemma \ref{brooksmatelskilemma}) and $H$ (large enough to ensure that each component of $\cone$ has length less than $\epsilon$), if the image of $\mathcal{I}_\alpha$ in $M_X(\alpha)$ intersected $\mathbf{T}_X^\epsilon(\alpha)$, then an essential closed curve in $\mathcal{I}_\alpha$ would be carried into $\mathbf{T}_X^{\epsilon_3}(\alpha)$, contrary to the fact that no power of $\gamma$ is conjugate into $\pi_1\big(\partial M(\alpha)\big)$.
So $\mathcal{I}_\alpha$ is disjoint from $\mathbf{B}_\alpha^\epsilon$.


Let $\mathcal{E}_\alpha$ be the component of 
\[
\overline{\QF\big(X, \sigma_{M(\alpha)}(X)\big) - \mathcal{C}_\alpha}
\]
 that faces $\sigma_{M(\alpha)}(X)$.
Let $E_t$ be the component of the boundary of the $t$--neighborhood of $\mathcal{C}_\alpha$ that lies in $\mathcal{E}_\alpha$.
The $E_t$ are convex and foliate $\mathcal{E}_\alpha$, see \cite{epsteinmarden}.
Let $Q_t$ be the submanifold of $\QF\big(X, \sigma_{M(\alpha)}(X)\big)$ facing $X$ whose boundary is $E_t$.

Since the map $F_{X,\alpha} \circ f_X$ is $L$--lipschitz and the diameter of $\F_X$ is bounded by a constant $A$ depending only on the diameter and thickness of $\sigma_M$, the diameter of $\F_{X,\alpha}$ is bounded by $LA$.
There is a constant $R$ depending only on $LA$ such that $\F_{X,\alpha}$ is contained in $Q_R$.
To see this, let $D$ be a universal lower bound on the diameter of all complete hyperbolic surfaces (as in the proof of Lemma \ref{oursurface}).
If $\F_{X,\alpha}$ lies outside of $Q_r$,  its diameter is at least that of $E_r$, which is at least $D \cosh r$.
Taking $r = \mathrm{arccosh}(2LA/D)$, we see that $\F_{X,\alpha}$ must intersect $Q_r$.
Taking $R= LA + \mathrm{arccosh}(2LA/D)$ guarantees that $\F_{X,\alpha}$ lies in $Q_R$.

Now, the image of $E_R$ in $M_X(\alpha)$ penetrates $\mathbf{T}_X^\epsilon(\alpha)$ to a depth no more than $R$.
By Theorem \ref{universal} and Brooks and Matelski's theorem \cite{brooksmatelski}, there is an $\hbar \geq H$ depending only on $R$ and $H$ such that the depth of $\mathbf{T}_X^\epsilon(\alpha)$ is at least $2R$ whenever $\ell(\alpha) \geq \hbar$.
So when $\ell(\alpha) \geq \hbar$, 
we may find a tubular neighborhood $\mathbf{T}_X^\omega(\alpha) \subset \mathbf{T}_X^\epsilon(\alpha)$ of $\cone$ whose preimage $\mathbf{B}_\alpha^\omega$ lies outside of $Q_R$.

Finally, let 
\[
N_\alpha^\omega \longrightarrow \QF\big(X, \sigma_{M(\alpha)}(X)\big) - \mathbf{B}_\alpha^\omega
\]
be the covering corresponding to 
\[
\pi_1(\mathcal{C}_\alpha) \
\subset \ \pi_1\Big(\QF\big(X, \sigma_{M(\alpha)}(X)\big) - \mathbf{B}_\alpha^\omega \Big)
\]
and let $N^\epsilon = \QF(X, \sigma_{M}(X)) -  \mathbf{B}^\epsilon$.
The map $F_{X,\alpha}$ lifts to an embedding 
\[
\widetilde F_{X,\alpha}\co N^\epsilon \to N_\alpha^\omega,
\]
and the composition 
\[
 \xymatrix{ \F_X \ar[r]^{\widetilde F_{X,\alpha}} & N_\alpha^\omega \ar[r] &  \QF\big(X, \sigma_{M(\alpha)}(X)\big) - \mathbf{B}_\alpha^\omega
}
\]
is our map $G_{X,\alpha}$.
As the tube $\mathbf{B}_\alpha^\omega$ misses it, the submanifold $Q_R \supset \F_{X,\alpha}$ lifts homeomorphically to a submanifold of $N_\alpha^\omega$, since $\pi_1(Q_R) = \pi_1(\mathcal{C}_\alpha)$.
Since $\F_X$ is embedded in $N^\epsilon$,  we conclude that $G_{X,\alpha}$  is an embedding when $\ell(\alpha) \geq \hbar$.
\end{proof}


Let $\hbar$ be as in the claim and assume that $\ell(\alpha) \geq \hbar$.

Since the principal curvatures of $\F_X$ are bounded away from zero by $1/2$, for sufficiently small $\eta$ the $\F_{X,\alpha}$ are strictly convex.
So there are homeomorphic normal projections $\Pi \co \sigma_M(X) \to \F_X$ and $\Pi_{\alpha} \co \sigma_{M(\alpha)}(X) \to \F_{X,\alpha}$. 
Since $\F_X$ and $\F_{X,\alpha}$ are smooth, these maps are diffeomorphisms.

We obtain a map $s_{X,\alpha} \co \sigma_M(X) \to \sigma_{M(\alpha)}(X)$ given by
\[
s_{X,\alpha} = \Pi_\alpha^{-1} \circ G_{X,\alpha} \circ \Pi  \, .
\]
As in the proof of Theorem \ref{collar}, we may write the derivative at a point $p$ as
\begin{align*}
Ds_{X,\alpha} & = A^{-1}B \\
&
=A^{-1}
\begin{pmatrix} \frac{1 + \kappa_1(\F_{X,\alpha})}{2}  &   0 \\
  0 &  \frac{1 + \kappa_2(\F_{X,\alpha})}{2}
\end{pmatrix}
A
\begin{pmatrix} w'  &   x' \\
  y' & z' 
\end{pmatrix}
\begin{pmatrix} \frac{2}{1 + \kappa_1(\F_X)}  &   0 \\
  0 &  \frac{2}{1 + \kappa_2(\F_X)}
\end{pmatrix} \\
& = A^{-1}
\begin{pmatrix} w\,  \frac{1 + \kappa_1(\F_{X, \alpha})}{1 + \kappa_1(\F_X)}  
&   x \,  \frac{1 + \kappa_1(\F_{X, \alpha})}{1 + \kappa_2(\F_X)} 
\\[7pt]
  y \,  \frac{1 + \kappa_2(\F_{X, \alpha})}{1 + \kappa_1(\F_X)} 
  & z \, \frac{1 + \kappa_2(\F_{X,\alpha})}{1 + \kappa_2(\F_X)} 
\end{pmatrix}
\end{align*}
where $A$ is the rotation carrying the principal directions of $\F_{X,\alpha}$ to $(1,0)$ and $(0,1)$---recall that we have chosen coordinates so that the principal directions of $\F_X$ are $(1,0)$ and $(0,1)$.
As before, we equip the tangent spaces to $\sigma_M(X)$ and $\sigma_{M(\alpha)}(X)$ at our chosen points with the standard Euclidean inner products.
Again, the numbers $w$, $x$, $y$, and $z$ are bounded in absolute value by a constant depending only on $L$.

There are now two cases depending on whether our point of interest in $\F_X$ is very nearly umbilic or not.

Let $\epsilon' > 0$.
There is a constant $\epsilon''> 0$ such that if the magnitude $| \, \kappa_1(\F_X) - \kappa_2(\F_X) \, |$ is less than $\epsilon''$ at $p$, and $\eta$ is small enough (so that the numbers $| \, \kappa_i(\F_X) - \kappa_i(\F_{X,\alpha})\, |$ are small), then the bilipschitz constant of $B$ is within $\epsilon'$ of that of $DG_{X,\alpha}$---for the simple reason that the entries of $B$ are then close to $w$, $x$, $y$, and $z$.

Suppose that $| \, \kappa_1(\F_X) - \kappa_2(\F_X) \, | \geq \epsilon''$.
Then, for any $\omega > 0$, we may choose $\eta$ very small in comparison to $\epsilon''$ so that the principal directions of $\F_{X,\alpha}$  are within $\omega$ of the images under $DG_{X,\alpha}$ of the principal directions of $\F_X$---this choice of $\eta$ depends only on $\omega$ and $\epsilon''$ by linearity of the second fundamental form and homogeneity of hyperbolic manifolds.
In this case $x$ and $y$ are very close to zero.
Since the $\kappa_i(\F_X)$ are both less than $2$, greater than $1/2$, and the numbers $| \, \kappa_i(\F_X) - \kappa_i(\F_{X,\alpha})\, |$ are small,
the bilipschitz constant of $B$ is within $\epsilon'$ of that of $DG_{X,\alpha}$.

In any case, the bilipschitz constant of $B$, and hence that of $Ds_{X,\alpha} = A^{-1}B$, is within $\epsilon'$ of that of $DG_{X,\alpha}$ provided $\ell(\alpha)$ is sufficiently large.
But the map $G_{X,\alpha}$ is bilipschitz with constant $L$.
We conclude that for $L$ very close to one and $\eta$ small, the Teichm\"uller distance between $\sigma_M(X)$ and $\sigma_{M(\alpha)}(X)$ is small in a manner independent of $X$.
\end{proof}

\section{Volume estimate}\label{volume}

In Section \ref{dehnfilling} we established the convention that given a manifold $M$ with conformal boundary $X$, a filling $M(\alpha)$ is assumed to be equipped with conformal boundary $X$.  
When working with manifolds with totally geodesic boundary, it is often useful to normalize so that the filled manifold has totally geodesic boundary. 
In the following theorem, a manifold with totally geodesic boundary is said to be obtained from a manifold $M$ with totally geodesic boundary via hyperbolic Dehn filling if it is topologically obtained from $M$ via Dehn filling and the core of the filling torus may be taken to be geodesic \textit{in the hyperbolic metric with totally geodesic boundary}.

\begin{theorem}[J\o rgensen]\label{jorg} Let $V > 0$. 
There is a finite list $N_1, \ldots, N_n$ of hyperbolic manifolds with closed totally geodesic boundary such that if $N$ is a hyperbolic manifold with totally geodesic boundary of volume less than $V$, then $N$ may be obtained from one of the $N_i$ by hyperbolic Dehn filling.
\end{theorem}
\begin{proof}[Sketch.] The theorem is typically stated for finite volume manifolds without boundary, as in Chapter 5 of \cite{thurston}, see also \cite{jorgmarden} and Chapter E.4 of \cite{bp}.

Let $N$ be a hyperbolic manifold of volume less than $V$ with totally geodesic boundary. 
Let $\mathcal{X}$ be a subset of the thick part $N^{\geq \epsilon}$ maximal with respect to the property that no two points of $\mathcal{X}$ have distance less than or equal to $\epsilon/2$.  
The bound on the volume bounds the cardinality of $\mathcal{X}$, as the $\epsilon/4$--neighborhood of $\mathcal{X}$ is embedded. 
Since $\mathcal{X}$ is maximal, the $\epsilon/2$--balls about the points of $\mathcal{X}$ cover $N^{\geq \epsilon}$.

Let $U$ be the universal cover of $N^{\geq \epsilon}$,  pull back our collection of $\epsilon/2$--balls to $U$, and consider the nerve $\widetilde{\mathcal{N}}$ of this covering of $U$---a convexity argument shows that the intersection of the elements of any subcollection of this covering is either empty or a topological cell.
  The group $\pi_1(N^{\geq \epsilon})$ acts freely on $\widetilde{\mathcal{N}}$, so the quotient $\mathcal{N} = \widetilde{\mathcal{N}}/\pi_1(N^{\geq \epsilon})$ is a simplicial complex and the quotient map $\widetilde{\mathcal{N}} \to \mathcal{N}$ is a covering.
  Since the balls in our original collection are bounded in number and size, there is a bound on the dimension and number of simplices in $\mathcal{N}$.
By a theorem of A. Weil, see pages 466--468 of \cite{corners}, the nerve $\widetilde{\mathcal{N}}$ is homotopy equivalent to $U$, which is contractible.
By Whitehead's Theorem, $\mathcal{N}$ is homotopy equivalent to $N^{\geq \epsilon}$.

We conclude that the number of fundamental groups of such $N^{\geq \epsilon}$ is finite.
As the $N^{\geq \epsilon}$ are acylindrical (see below), they are determined up to homeomorphism by their fundamental groups, thanks to a theorem of K. Johannson \cite{johannson}.
It follows that $N^{\geq \epsilon}$ belongs to a finite list of manifolds depending only on $V$.

If a Margulis tube in the thin part $N^{\leq \epsilon}$ is incident to $\partial N$, there is only one topological Dehn filling there, and as the geometry of $N$ depends only on its topology, by Mostow--Prasad Rigidity \cite{prasad}, we are free to fill these tubes to obtain a manifold $N^*$. 

To see that this manifold admits a hyperbolic structure with totally geodesic boundary, note that the double $DN$ is a hyperbolic manifold of finite volume.  
Now, the complement of a collection of simple geodesics in a finite volume hyperbolic $3$--manifold admits a finite volume hyperbolic metric itself, see \cite{kojimanonsingular}, for instance, and so $DN^*$ admits such a metric. 
The involution on $DN^*$ implies that $\partial N^*$ is totally geodesic in this metric, and we obtain a metric with totally geodesic boundary on $N^*$.
\end{proof}

\begin{theorem}[Volume Estimate]\label{estimate} Let $M$ be a finite volume hyperbolic $3$--manifold  with nonempty closed totally geodesic boundary. 
Then there are $A, B, \delta > 0$ depending only on the volume of $M$ such that $\sigma_M\big(\T(\partial M)\big)$ is $\delta$--thick and 
$
B \leq \diam(\sigma_M) \leq A.
$
\end{theorem}
\begin{proof} Given a finite set $\mathfrak{M}$ of finite volume manifolds, let $\mathfrak{C} (\mathfrak{M})$ be the maximum number of cusps possessed by any element of $\mathfrak{M}$.

 Let $V >0$. 

Let $\mathfrak{M}_1 = \mathfrak{M}_1(V) = \{ M_1, \ldots, M_{n_1} \}$ be the list of manifolds given by J\o rgensen's theorem after discarding any manifolds with empty boundary.

For each $i$, let $\delta_i$, $A_i$, and $B_i$ be positive constants such that the image of $\sigma_{M_i}$ is $\delta_i$--thick and 
\[
B_i \leq \diam(\sigma_{M_i}) \leq A_i ,
\]
as we may by Theorem \ref{nonconstant} and the Bounded Image Theorem.
Let $\epsilon_i > 0$ be small enough so that the $\epsilon_i$--neighborhood of the $\delta_i$--thick part of $\T(\overline{\partial M_i})$ is $\delta_i/2$--thick,
\[
\frac{B_i}{2} \leq B_i - \epsilon_i \, ,
\]
and
\[
A_i + \epsilon_i \leq 2A_i.
\]

For each $i$, let $\mathcal{E}_i$ be the set of  filling slopes 
$\alpha$ for $M_i$ such that the normalized length of \textit{every} component of $\alpha$ is less than the constant $\hbar_{\epsilon_i}(M_i)$ given by the Filling Theorem \textit{and} 
such that there is a hyperbolic Dehn filling $M_i(\alpha)(\beta)$ of $M_i(\alpha)$ with totally geodesic boundary.
Being the complement of a union of geodesics in $M_i(\alpha)(\beta)$, the manifold $M_i(\alpha)$ admits a hyperbolic structure with totally geodesic boundary---again, see \cite{kojimanonsingular}.

For each $\alpha$ in $\mathcal{E}_i$, adjoin the manifold $M_i(\alpha)$ to the list $\mathfrak{M}_1$ to obtain a finite list of manifolds $\mathfrak{M}_2 = \{M_1, \ldots, M_{n_1}, M_{n_1+1}, \ldots, M_{n_2}\}$. 
Note that $\mathfrak{C} (\mathfrak{M}_2 - \mathfrak{M}_1) < \mathfrak{C}(\mathfrak{M}_1)$.
For $n_1 +1 \leq i \leq n_2$, we define $\epsilon_i$, as before, to be small enough so that the $\epsilon_i$--neighborhood of the $\delta_i$--thick part of $\T(\overline{\partial M_i})$ is $\delta_i/2$--thick,
$
B_i/2 \leq B_i - \epsilon_i \, ,
$
and
$
A_i + \epsilon_i \leq 2A_i ,
$
where the image of $\sigma_{M_i}$ is $\delta_i$--thick and 
$
B_i \leq \diam(\sigma_{M_i}) \leq A_i .
$

We apply the argument again to the set $\mathfrak{M}_2$ to obtain a set $\mathfrak{M}_3$.
Note that any element in $\mathfrak{M}_3 - \mathfrak{M}_2$ is obtained by filling a manifold in $\mathfrak{M}_2 - \mathfrak{M}_1$, and so 
\[
\mathfrak{C}(\mathfrak{M}_3 - \mathfrak{M}_2) < \mathfrak{C}(\mathfrak{M}_2 - \mathfrak{M}_1).
\]
Continuing in this manner, we obtain a sequence of finite sets 
\[
\mathfrak{M}_1 \subset \mathfrak{M}_2 \subset \cdots \subset \mathfrak{M}_i \subset \cdots
\]
with 
\[
\mathfrak{C}(\mathfrak{M}_{i} - \mathfrak{M}_{i-1}) > \mathfrak{C}(\mathfrak{M}_{i+1} - \mathfrak{M}_i)
\]
for all $i$.
So this process terminates in a finite set $\mathfrak{M}_k = \{M_1, \ldots, M_{n_k} \}$.

Now, by construction, this set has the property that given a manifold $M$ with totally geodesic boundary and volume less than $V$, there exists an 
$i$ in $\{ 1, \ldots, n_k\}$ such that $M = M_i(\gamma)$ where the cores of the filling tori are geodesics in $M$ and the normalized length of each component of $\gamma$ is at least $\hbar_{\epsilon_i}(M_i)$.
So, letting
\[
\delta = \min\bigg\{ \frac{\delta_i}{2} \ \bigg| \ 1 \leq i \leq n_k \bigg\} ,
\]
\[
B = \min\bigg\{ \frac{B_i}{2} \ \bigg| \ 1 \leq i \leq n_k \bigg\} ,
\]
and
\[
A = \max \{2A_i \ | \ 1 \leq i \leq n_k \},
\]
we see that for any $M$ with totally geodesic boundary and volume no more than $V$, the image of $\sigma_M$ is $\delta$--thick and 
$
B \leq \diam(\sigma_M) \leq A.
$
\end{proof}

\section{Pinching}\label{pinch}

By Fujii and Soma's Density Theorem \cite{fujiisoma}, there are hyperbolic $3$--manifolds whose totally geodesic boundaries contain pants decompositions as short as desired, see also \cite{kentboundaries}, and so the following theorem produces manifolds whose skinning maps have diameters descending to zero.

\begin{theorem}[Bromberg--Kent]\label{brombergkent} Let $S$ be a closed hyperbolic surface.
For each $\epsilon > 0$ there is a $\delta > 0$ such that if $M$ is a hyperbolic $3$--manifold with totally geodesic boundary $\Sigma$ homeomorphic to $S$ containing a pants decomposition $\mathcal{P}$ each component of which has length less than $\delta$, then the diameter of $\sigma_M$ is less than $\epsilon$.
\end{theorem}
\begin{proof} For simplicity we assume that $S$ is connected of genus $g$.
We let $\Gamma$ denote the unifomizing Kleinian group for our $3$--manifold, and as we work with the open manifold $\HH^3/\Gamma$, we let $M = \HH^3/\Gamma$ and consider $\Sigma$ a totally geodesic surface inside $M$. 

Let $\mathbf{T}_\Sigma$ be the $\epsilon_3$--Margulis tube about $\mathcal{P}$ in $M$.

Now consider the double $D\mathcal{C}(M)$ of the convex core of $M$. 
The tube $\mathbf{T}_\Sigma$ isometrically embeds in the double and so we call its image by the same name.

Let $m = m_1 \sqcup  \cdots \sqcup m_{3g-3}$ be a geodesic $1$--manifold in $\partial \mathbf{T}_\Sigma$  (equipped with the induced path metric) such that each component $m_i$ bounds a disk in $\mathbf{T}_\Sigma$. 
Each component $m_i$ of $m$ is a \textbf{meridian}. 
If the total length of $\mathcal{P}$ is small enough, then the depth of each component of $\mathbf{T}_\Sigma$ is large, by Brooks and Matelski's Theorem \cite{brooksmatelski}. 
It follows that the \textit{actual} lengths of the meridians $m_i$ in $\partial \mathbf{T}_\Sigma$ are all large---this is true of the \textit{meridian} of any deep Margulis tube about a geodesic.

The intersection of the surface $\partial \mathcal{C}(M)$ and $\mathbf{T}_\Sigma$ defines a slope in $\partial \mathbf{T}_\Sigma$ whose components we call the \textbf{longitudes}.
Notice that the actual lengths of the longitudes are bounded above as $\epsilon$ tends to zero, as they are trapped in geodesic pairs of pants whose boundary components have length $\epsilon$.

By Brock and Bromberg's Drilling Theorem, Theorem \ref{drilling} here, we may drill $\mathcal{P}$ out from $D\mathcal{C}(M)$ to obtain a complete hyperbolic manifold $M'$  and an $L$--bilipschitz map 
\[
f \co D\mathcal{C}(M) - \mathbf{T}_\Sigma \longrightarrow M' - \mathbf{P}
\]
where $L$ is very close to one and $\mathbf{P}$ is a Margulis tube about the $3g-3$ new rank--two cusps in $M'$---here we are using the Drilling Theorem in the case when the conformal boundary is empty.
Letting $m_1', \ldots , m_{3g-3}'$ denote the geodesic representatives of the $f(m_1), \ldots , f(m_{3g-3})$, respectively, we see that the actual flat lengths of the $m_i'$ are large in $\partial M'$. 
Also note that the longitudes are not much longer in $\partial \mathbf{P}$ than they were in $\mathbf{T}_\Sigma$.

The manifold $M'$ is the double of a manifold whose interior is homeomorphic to $M$ and whose boundary is a union $\mathcal{S} = S_1 \cup \cdots \cup S_{2g-2}$ of thrice--punctured spheres. 
The involution on $M'$ implies that $\mathcal{S}$ may be isotoped to be totally geodesic in $M'$. 
So the surface $\mathcal{S}$ cuts the union of tori $\partial \mathbf{P}$ into a union $\mathbf{A}$ of annuli.
It also cuts each meridian $m_i'$ into two arcs of equal length.
If $\zeta$ is one of these arcs and $A$ its annulus in $\mathbf{A}$, then the involution implies that $\zeta$ intersects $\partial A$ in right angles.
Since the longitudes are bounded above in length, and $\zeta$ is long, it follows that the conformal modulus of $A$ is large.

A separate application of the Drilling Theorem allows the drilling of $\mathcal{P}$ from $M$ to obtain a complete hyperbolic manifold $N$ with conformal boundary $\Sigma$ and an $L$--bilipschitz map 
\[
f \co M - \mathbf{T}_\Sigma \longrightarrow N - \mathbf{P}_\Sigma
\]
where $L$ is very close to one and $\mathbf{P}_\Sigma$ is a Margulis tube about the new rank--two cusps. 

The manifold $N$ contains $2g-2$ incompressible properly embedded thrice--punc-tured spheres $S_1, \ldots, S_{2g-2}$ whose ends lie in $\mathbf{P}_\Sigma$.  
A theorem of Adams \cite{adams} says that we may take the $S_i$ to be totally geodesic in \textit{any} complete hyperbolic structure on $N$---in \cite{adams}, the theorem is stated for finite volume hyperbolic manifolds, but it is easily seen that the proof is valid as long as the ends of the punctured spheres lie in the cusps of the ambient manifold.

Given an Riemann surface $X$ in $\T(S) \cong \GF(N)$, let $\mathcal{S}_X$ be the totally geodesic representative of $S_1 \cup \cdots \cup S_{2g-2}$ in $N_X$ given by Adams' theorem.
Let $N_X'$ be the closure of the component of $N_X - \mathcal{S}_X$ homeomorphic to $M$, and let $N'= N_\Sigma'$.
Since thrice--punctured spheres have no moduli, all of the $N_X'$ are isometric to $N'$, by Theorem \ref{geomfinite}, and so we have isometric embeddings
\[
\iota_X \co N' \to N_X.
\]

Note that our manifold $M'$ above is the double of $N'$.

For $X$ in $\T(S)$, let $\mathbf{P}_X$ be a union of cuspidal $\epsilon_3$--Margulis tubes in $N_X$.
Since $N'$ is isometrically embedded in $N_X$, the intersection $\partial \mathbf{P}_X \cap N'$ is a collection $\mathbf{B}$ of horospherical annuli perpendicular to $\mathcal{S}$.
Moreover, each annulus in $\mathbf{B}$ is conformally equivalent to the corresponding annulus in $\mathbf{A}$.

Being trapped in finite volume geodesic pairs of pants, the longitudes in $\partial \mathbf{P}_X$ have lengths bounded above independent of $X$.
Since the Margulis constant $\epsilon_3$ is universal, their lengths are uniformly bounded below as well.

Now, if $n(X)$ is a geodesic meridian on $\partial \mathbf{P}_X$, the surface $\mathcal{S}$ cuts it into two arcs.  
Let $\zeta$ be the arc contained in $N'$, and let $B$ be the annulus in $\mathbf{B}$ containing it.
Since $B$ and the corresponding annulus $A$ in $\mathbf{A}$ are conformally equivalent, they have the same conformal modulus.
Since the modulus is large, and the actual length of the longitude is bounded \textit{below}, the actual length of $\zeta$ must be large, and it follows that the actual length of $n(X)$ is large as well.
Since the actual lengths of the longitudes are bounded above, we conclude that the normalized length of $n(X)$ is large for all $X$.

If the normalized lengths of these meridians are large enough (which is guaranteed if the total hyperbolic length of $\mathcal{P}$ is small enough), we may perform hyperbolic Dehn filling at all of them to obtain $M_X$. 
Moreover, the Drilling Theorem again tells us that for each $X$ there is an $L$--bilipschitz map
\[
g_X \co N_X - \mathbf{P}_X \longrightarrow M_X - \mathbf{T}_X  .
\]

\bigskip
\noindent
We now construct a surface $\mathfrak{G}$ in $N'$ with diffeomorphic Gauss map that will allow us to compare skinning surfaces.  
Unlike the surfaces in the proofs of Theorems \ref{collar} and \ref{fill}, this surface will not be convex.
The idea is as follows. 
We would like to construct a surface by cutting the surface $\mathcal{S}$ along $\partial \mathbf{P}_\Sigma$, tubing the resulting boundary components together by annuli in $\partial \mathbf{P}_\Sigma$, and smoothing.
Unfortunately, the hyperbolic Gauss map of such a surface will not be a diffeomorphism, as it may have principal curvatures below $-1$.
To correct this, we carefully extend the surface $\mathcal{S} - \mathbf{P}_\Sigma$ into the cusps $\mathbf{P}_\Sigma$ in a ``rotationally" symmetric way \textit{before} tubing to get a surface with diffeomorphic Gauss map after smoothing.
\begin{figure}
\begin{center}
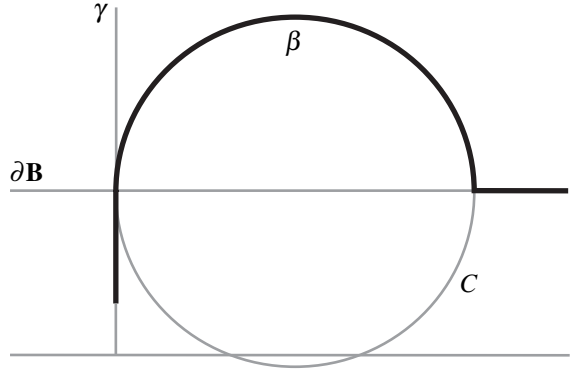
\end{center}
\caption{Building a slice of $\mathfrak{G}$ near $\mathbf{P}_\Sigma$.}
\label{nonconvexsurface}
\end{figure}

Outside a small neighborhood $U$ of $\mathbf{P}_\Sigma$, the surface $\mathfrak{G}$ is simply $\mathcal{S} - U$. 
To describe $\mathfrak{G}$ in $U$, we begin by defining a meridional slice of $\mathfrak{G}$.

Pick a component $\mathbf{P}_0$ of $\mathbf{P}_\Sigma$ and consider one of its preimages in $\HH^3$ normalized to be the horoball $\mathbf{B}$ centered at infinity that passes through $i$. 
Choose a sphere $S_i \subset \mathcal{S}$ that intersects $\mathbf{P}_0$
and further normalize so that a preimage $\widetilde S_i$ of $S_i$ contains the geodesic $\gamma$ passing through zero and infinity.
There is a hyperbolic plane $H$ that intersects both $\partial \mathbf{B}$ and $\widetilde S_i$ at right angles.

The plane $H$ is pictured in Figure \ref{nonconvexsurface} with $\partial \mathbf{B}$ and $\gamma$ in gray. 
Also in gray is a circle $C$ tangent to $\gamma$ at $i$ and whose Euclidean radius is slightly larger than one---note that, choosing inward pointing normals, the hyperbolic curvature of $C$ is strictly greater than $-1$. 
The black arc $\beta$ is the concatenation of an arc in $\gamma = \widetilde S_i \cap H$, an arc in $C$, and an arc of $\partial \mathbf{B}$.
Note that if the meridian in $N$ is long enough, then $\beta$ may be taken disjoint from all of the other preimages of components of $\mathcal{S}$.
If the meridian is still longer, then two adjacent lifts of components of $\mathcal{S}$ that intersect $\mathbf{B}$ may be joined be an embedded arc formed by concatenating a copy of $\beta$ and a copy of its reflection through a vertical line---this will guarantee that our surface $\mathfrak{G}$ is embedded---see Figure \ref{nonconvexschematic}.

Now, the arc $\beta$ is not smooth, but a small perturbation changes $\beta$ into a smooth arc $\alpha$ that agrees with $\beta$ near its endpoints and all of whose curvatures are strictly greater than $-1$.  
Now, pushing $\alpha$ down into our neighborhood $U$ of $\mathbf{P}_0$, revolving it through the product structure there, and taking the union with the surface $\mathcal{S} - U$, we obtain a smooth embedded surface $\mathfrak{G}$ all of whose normal curvatures are strictly greater than $-1$. 
See Figure \ref{nonconvexschematic} for a schematic of $\mathfrak{G}$ in $M$.
\begin{figure}
\begin{center}
\includegraphics{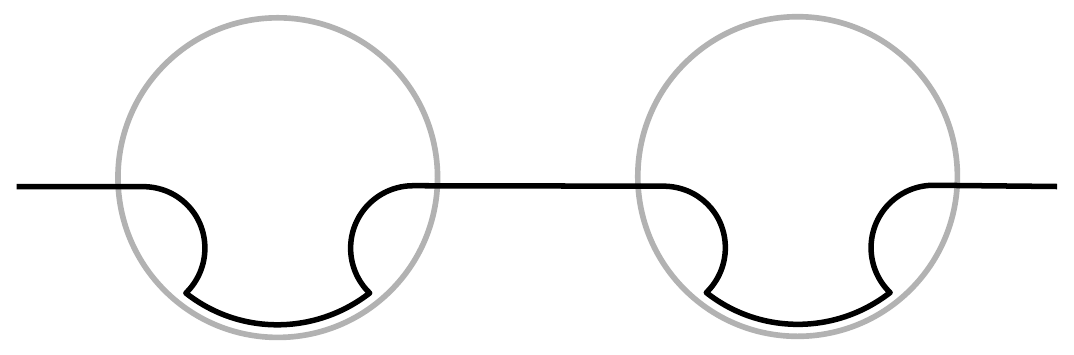}
\end{center}
\caption{The surface $\mathfrak{G}$ intersecting two Margulis tubes.
Here, the skinning surface lies ``below" $\mathfrak{G}$.}
\label{nonconvexschematic}
\end{figure}

By construction, $\mathfrak{G}$ is contained in the $\epsilon'$--thick part of $N$ for some universal $\epsilon'$. 
Note that the surfaces $f^{-1}(\mathfrak{G})$ and $g_X \circ \iota_X (\mathfrak{G})$ are homotopic to the convex core boundaries of $M$ and $M_X$, respectively.
Theorem \ref{universal} and Proposition \ref{curvdrilling} tell us that the normal curvatures of  $f^{-1}(\mathfrak{G})$ and $g_X \circ \iota_X(\mathfrak{G})$ are very close to those of $\mathfrak{G}$, and, in particular, are greater than $-1$. 
This is all that is needed for the existence of diffeomorphic normal projections $\Pi \co \overline{\Sigma} \to f^{-1}(\mathfrak{G})$ and $\Pi_X \co \sigma_M(X) \to g_X  \circ \iota_X(\mathfrak{G})$, as the hyperbolic Gauss map of a smooth embedded surface in $\HH^3$ is a diffeomorphism provided all of its normal curvatures are strictly greater than $-1$. 

As in the proof of Theorem \ref{fill}, as we let the normal curvatures of $g_X  \circ \iota_X (\mathfrak{G})$ and $f^{-1}(\mathfrak{G})$ tend to those of $\mathfrak{G}$, letting our bilipschitz constant $L$ tend to one, we see that the composition
\[
\Pi_X^{-1} \circ g_X  \circ \iota_X \circ f \circ \Pi \co \overline{\Sigma} \to \sigma_M(X)
\]
is very close to conformal in a manner independent of $X$.
\end{proof}

\section{No universal bound}\label{nouniversal}

\begin{theorem}[Flexible Manifolds]\label{flex} For any $R >0$, 
there is a small $d$ such that if $M$ is a hyperbolic $3$--manifold with connected totally geodesic boundary $\Sigma$ such that every metric collar about $\Sigma$ has depth at most $d$, then the 
skinning map of $M$ has diameter greater than $R$.
\end{theorem}
\begin{proof}  Let $N_n$ be any sequence of hyperbolic $3$--manifolds with connected totally geodesic boundary $\Sigma_n$ such that every collar about $\Sigma_n$ has depth at most $d_n$ with $\lim \, d_n = 0$.
Let $\Gamma_n$ be a uniformizing Kleinian group for $N_n$, and let $\sigma_n = \sigma_{N_n}$. 
We will show that the diameters of the $\sigma_n$ tend to infinity.

After normalizing, we may assume that $\U$ is the component of the domain of discontinuity of $\Gamma_n$ that uniformizes $\Gamma_{\Sigma_n}$. 
Since $\lim \, d_n = 0$, we may further normalize so that there are loxodromic elements $g_n$ in the $\Gamma_n$ with the property that
\begin{equation}\label{areadecay}
\lim_{n \to \smallinfinity} \mathrm{meas}(\ELL - g_n \U) = 0  .
\end{equation}
This normalization produces the picture in Figure \ref{domain}. 
There is a fundamental domain for $\Gamma_{\overline{\Sigma}_n}$ in $\ELL$ that contains $g_n \U$, and so, by \eqref{areadecay},
the $\Gamma_{\Sigma_n}$ have fundamental domains $\omega_n$ in $\U$ such that
\[
\lim_{n \to \smallinfinity} \mathrm{meas}(g_n \U - g_n \omega_n) = 0.
\]
With \eqref{areadecay}, we have
\begin{equation}\label{fundamentaldomains}
\lim_{n \to \smallinfinity} \mathrm{meas}(\ELL - g_n \omega_n) = 0.
\end{equation}

Let $k$ lie in the interval $(0,1)$ and let $K = (1+k)/(1-k)$.
Define a Beltrami differential $\mu_n$ for $\Gamma_{\Sigma_n}$ by declaring that
\[
\mu_n(z) = k \, \frac{\overline{g_n'(z)}}{g_n'(z)} 
\]
on $\omega_n$ and extending to $\U$ equivariantly.

The constant function $1/\pi$ is a quadratic differential on $\ELL$.
So
\begin{align}
\frac{1}{\mathrm{K}\big[ \sigma_{n}([\mu_{n}])\big]_{\ELL}} 
& \leq \iint_{\ELL} \frac{1}{\pi} \frac{\ |1 -  \sigma_{n}(\mu_{n}) |^2}{1 - | \sigma_{n}(\mu_{n})|^2}  \ \dee x \, \dee y  \quad \quad \ \  \mathrm{by \ \eqref{firstineq}} \notag \\
& = \iint_{g_n  \omega_n}  \frac{1}{\pi} \frac{\ |1 -  \sigma_{n}(\mu_{n})|^2}{1 - | \sigma_{n}(\mu_{n})|^2}   \ \dee x \, \dee y  \notag \\
& \quad \quad \quad \quad \quad \quad \quad \quad \quad +  \iint_{\ELL - g_n \omega_n} \frac{1}{\pi} \frac{\ |1 -  \sigma_{n}(\mu_{n})|^2}{1 - | \sigma_{n}(\mu_{n})|^2}   \ \dee x \, \dee y  \notag \\
& \leq  \iint_{\omega_n}  \frac{1}{\pi}
\frac{\
\Big | 1 -   k \frac{\, \overline{g_n'} \, g_n' \, }{ g_n' \, \overline{g_n'}} 
\Big |^2
}{1 -  k^2}   |g_n'|^2 \ \dee x \, \dee y  
 \notag \\
& \quad \quad \quad \quad \quad \quad \quad \quad \quad + \frac{{1 + k}^2}{\pi(1 - k^2)} \, \mathrm{meas}(\ELL - g_n \omega_n) \notag \\
& = \frac{1}{K} \iint_{\omega_n}  \frac{1}{\pi} \, |g_n'|^2
   \ \dee x \, \dee y  
 \notag \\
&   \quad \quad \quad \quad \quad \quad \quad \quad \quad + \frac{{1 + k}^2}{\pi(1 - k^2)} \, \mathrm{meas}(\ELL - g_n \omega_n)  \label{secondterm}
\end{align}

\noindent By \eqref{fundamentaldomains}, the term \eqref{secondterm} is less than $\epsilon$ for all sufficiently large $n$, and so
\[
\frac{1}{\mathrm{K}\big[ \sigma_{n}([\mu_{n}])\big]_{\ELL}}  \leq \frac{1}{K} + \epsilon ,
\]  
since 
\[
 \iint_{\U}  \frac{1}{\pi} \, |g_n'|^2
   \ \dee x \, \dee y  \leq 1.
\]

\noindent Putting $k$ close to one---so that $K$ is large---we see that we may take 
\[
\mathrm{K}\big[ \sigma_{n}([\mu_{n}])\big]_{\ELL} \leq \mathrm{K}\big[ \sigma_{n}([\mu_{n}])\big]_{\overline{\Sigma}_n}
\]
 as large as we like.   
 This means that the distance 
 $\dee_{\T(\overline{\Sigma}_n)}\big([0], \sigma_{n}([\mu_{n}])\big)$ is large, and so the skinning map $\sigma_{n}$ has large diameter. 
\end{proof}

\section{The Bounded Image Theorem}\label{bddimage}

In the following $M$ will denote a compact oriented irreducible atoroidal $3$--manifold with incompressible boundary of negative Euler characteristic. 
Let $\partial_0 M$ be the part of the boundary that contains no tori. 
\begin{theorem}[Thurston]\label{continuous} If $M$ is acylindrical, the skinning map $\sigma_M$ admits a continuous extension
\[
\sigma_M \co \AH(M) \to \T(\overline{\partial_0 M}) .
\]
\end{theorem}
\noindent Since $\AH(M)$ is compact, this implies the Bounded Image Theorem. 

\bigskip
\noindent To complete the final gluing step in the Geometrization of Haken Manifolds, Thurston originally suggested the following strong form of the Bounded Image Theorem---see statement ($\mathrm{II'}$) on page 85 of \cite{morgansurvey}. 
\begin{statement} Let $M$ be an orientable irreducible $3$--manifold with incompressible boundary that is not an interval bundle over a surface. 
Let $\tau \co \overline{\partial_0 M} \to\partial_0 M$ be a homeomorphism such that the $3$--manifold $M/\tau$ is atoroidal. 
Then there is a number $n$ such that the map $(\tau \circ\sigma_M)^n \co \T(\partial_0 M) \to \T(\partial_0 M)$ has bounded image.
\end{statement}
\noindent We know of no proof of this global statement. 
In all accounts of Thurston's argument, a local statement is used at the final gluing step, namely that for any $X$ in $\T(\partial_0 M)$, the sequence $(\tau \circ\sigma_M)^n(X)$ is bounded. 
See \cite{morgansurvey,skinmcmullen,otalhaken,kapovichbook}.

\bigskip
\noindent We continue to the proof of Theorem \ref{continuous}. 

Let $N$ be a manifold.  
A \textbf{compact core} of $N$ is a compact submanifold with the property that the inclusion is a homotopy equivalence.  
It is a theorem due to G. P. Scott and P. Shalen that aspherical $3$--manifolds with finitely generated fundamental groups have compact cores \cite{scott}.

If $N$ is a hyperbolic $3$--manifold, let $N^0$ denote the complement of its $\epsilon_3$--cuspidal thin part. 
An end $E$ of $N^0$ with a neighborhood $U$ homeomorphic to $S \times [0,\infinity)$ is \textbf{simply degenerate}   if there is a sequence of pleated surfaces $f_n \co \Sigma_n \to U$ homotopic in $U$ to the inclusion $S \times \{0 \} \to U$ whose images leave every compact set in $U$.
As is common, we will often refer to a neighborhood of an end \textit{as} the end.

To show that the image of our map lies in the Teichm\"uller space, we will need Thurston's Covering Theorem, Theorem 9.2.2 of \cite{thurston}, see \cite{canary}:

\begin{covering} Let $f \co N \to M$ be a locally isometric covering of hyperbolic $3$--manifolds such that $\pi_1(N)$ is finitely generated. 
Let $E$ be  a simply degenerate end of $N^0$ corresponding to an incompressible surface in $N^0$. Then either
\begin{enumerate}
\item $E$ has a neighborhood $U$ such that $f$ is finite-to-one on $U$, or
\item $M$ has finite volume and has a finite cover $M'$ that fibers over the circle and $N$ is finitely covered by the cover of $M'$ corresponding to the fiber subgroup. \qed
\end{enumerate} 
\end{covering}
\noindent Canary's strong version of this theorem removes the requirement that the end correspond to an incompressible surface, provided the manifolds in question are tame.  
Tameness is provided in Thurston's theorem by a theorem of Bonahon \cite{bonahon}, and is now known for all hyperbolic manifolds with finitely generated fundamental group by the recent solution to Marden's Tameness Conjecture by I. Agol \cite{agol} and D. Calegari and D. Gabai \cite{calgabai}.

For the continuity we will need the following embedding theorem of Canary and J. Anderson, which is implicit in \cite{andersoncanary}.

\begin{theorem}[Anderson--Canary \cite{andersoncanary}]\label{embeddedcore} Let $M$ be a hyperbolic $3$--manifold with finitely generated fundamental group and totally geodesic boundary. 
Let $M_n \in \AH(M)$ be a sequence of manifolds converging algebraically to $M_\smallinfinity^a$ and geometrically to $M_\smallinfinity^g$. 
Then there is a compact core $\K$ of $M_\smallinfinity^a$ such that the restriction of the covering $M_\smallinfinity^a \to M_\smallinfinity^g$ to $\K$ is an embedding. 
\end{theorem}
\begin{proof} Let $\Gamma^a$ and $\Gamma^g$ be the Kleinian groups for $M_\smallinfinity^a$ and $M_\smallinfinity^g$, respectively.
By tameness and Proposition 2.7 of \cite{ACCS}, for any element $\gamma$ of $\Gamma^g - \Gamma^a$, the intersection $\Lambda_{\Gamma^a} \cap \gamma \Lambda_{\Gamma^a}$ has cardinality at most one.
So either $\Gamma^g = \Gamma^a$ or $\Gamma^a$ has a nonempty domain of discontinuity.
In the first case, the theorem is obvious.  
In the second, the theorem follows from Corollary B of \cite{andersoncanary}.
\end{proof}

\begin{proof}[Proof of Theorem \ref{continuous} (Brock--Kent--Minsky)]
For simplicity we assume that only one component $\partial_0 M$ of $\partial M$ has negative Euler characteristic.

\bigskip \noindent \textit{The map.}
Let $M_\rho$ be a hyperbolic manifold in $\AH(M)$,  with end invariant $\lambda$, corresponding to the character of a representation $\rho \co \pi_1(M) \to \SL_2 \C$. 
Let $f \co N_\rho \to M_\rho$ be the cover of $M_\rho$ corresponding to $\partial_0 M$.  
This manifold has two end invariants, one of them $\lambda$ corresponding to the lift of the end of $M_\rho$, and a skinning invariant $\sigma_M(\rho)$ corresponding to the new end $E$ (the end $E$ is tame by Bonahon's theorem).

We claim that $E$ has no rank--one cusp.
To see this, suppose to the contrary that there is a rank--one cusp $\mathbf P$ in $E$.
If $\mathbf P$ covers a rank--two cusp in $M_\rho$, we conclude that a closed curve in $\partial_0 M$ is homotopic into a torus in $\partial M$,
contrary to the fact that $M$ is acylindrical. 
So $\mathbf P$ must cover a rank--one cusp in $M_\rho$.
But then, by considering the corresponding annulus in $N_\rho$ and its image in $M_\rho$, we conclude that there is a nontrivial conjugacy in $\pi_1(M)$ between two elements of $\pi_1(\partial_0 M)$, again contradicting the fact that $M$ is acylindrical. 
So $E$ is an end of $N_ \rho^0$. 
Furthermore, it is either simply degenerate or conformally compact.

Now, the restriction of $f$ to any neighborhood of $E$ is infinite-to-one, and so, by the Covering Theorem, $E$ cannot be simply degenerate, as $M_\rho$ has infinite volume.  
This means that $\sigma_M(\rho)$ is a Riemann surface homeomorphic to $\partial_0 M$.
Now, it may happen that the orientation of $M_\rho$ naturally places its skinning invariant in $\T(\partial_0 M)$ rather than $\T(\overline{\partial_0 M})$, in which case we let $\sigma_M(\rho)$ be the mirror image of this invariant.
We thus obtain a map
\[
\sigma_M \co \AH(M) \to \T(\overline{\partial_0 M})
\]
extending the skinning map, defined by $\sigma_M(M_\rho) = \sigma_M(\rho)$.

\bigskip \noindent \textit{Continuity.}
To see that $\sigma_M$ is continuous, let $M_n$ be a sequence in $\AH(M)$ converging algebraically to $M_\smallinfinity^a$.  
Let $N_n$ be the cover of $M_n$ corresponding to $\partial_0 M$.
Since the map 
\[
\AH(M) \to \AH(\partial_0 M)
\]
induced by inclusion is continuous (being the restriction of the regular function induced on $\SL_2 \C$--character varieties), the $N_n$ converge in $\AH(\partial_0 M)$ to a manifold $N_\smallinfinity^a$.

We pass to a subsequence so that the $M_n$ converge geometrically---see Corollary 9.1.8 of \cite{thurston} and Proposition 3.8 of \cite{jorgmarden2}.  
After passing to a deeper subsequence, this produces geometric convergence of the $N_n$, and considering that this is convergence of the uniformizing Kleinian groups in the Chabauty topology, we have the following commutative diagram of covering spaces:
\[
 \xymatrix{N_\smallinfinity^a \ar[d]_{}  \ar[dr]^{} &   \\
  M_\smallinfinity^a  \ar[dr] & N_\smallinfinity^g \ar[d] \\
  & M_\smallinfinity^g
  } 
\]

Let $\Gamma_n$, $\Gamma_\smallinfinity^a$, and $\Gamma_\smallinfinity^g$ be the uniformizing Kleinian groups for $N_n$, $N_\smallinfinity^a$, and $N_\smallinfinity^g$, respectively.

By Theorem \ref{embeddedcore}, there is a core $\K^a$ of $M_\smallinfinity^a$ that embeds in $M_\smallinfinity^g$. 
We call its image $\K^g$.
We let $\K_n$ be the image of the core $\K^g$ in the approximate $M_n$.

Now, a peripheral $\pi_1$--injective surface $\mathcal{S} \to M$ in an acylindrical manifold $M$ has a unique non-simply connected lift $\widetilde{\mathcal{S}} \to \widetilde{M}$ to the cover corresponding to the boundary, and we let $\mathfrak{S}^a$ denote this \textit{bona fide} lift of $\partial \K^a$ to $N_\smallinfinity^a$.
We let $\mathfrak{S}^g$ denote the image of $\mathfrak{S}^a$ in $N_\smallinfinity^g$.
Diagrammatically:
\[
 \xymatrix{\mathfrak{S}^a \ar[d]_{}  \ar[dr]^{\cong} &   \\
  \K^a  \ar[dr]_{\cong} & \mathfrak{S}^g \ar[d] & \ar[r] & \\
  & \K^g &
  } 
\quad \quad \quad
  \xymatrix{N_\smallinfinity^a \ar[d]_{}  \ar[dr]^{} &   \\
  M_\smallinfinity^a  \ar[dr] & N_\smallinfinity^g \ar[d] \\
  & M_\smallinfinity^g
  } 
\]
We let $\mathfrak{S}_n$ be the image of $\mathfrak{S}^g$ in $N_n$.

The algebraic limit $N_\smallinfinity^a$ has two ends: one the isometric lift of the end of $M_\smallinfinity^a$, the ``left" side; and another we call $E^a$, the ``right" side---for psychological reasons, we are choosing a homeomorphism $N_\smallinfinity^a \cong \R \times \partial_0 M$.
Since $\sigma_M(M_\smallinfinity^a)$ lies in $\T(\overline{\partial_0 M})$, the end $E^a$ is conformally compact.

Each $\mathfrak{S}_n$ cuts each $N_n$ into two pieces $A_n$  and $E_n$: the first facing the end lifted from $M_n$; the second facing the skinning surface of $M_n$. 
The geometric convergence implies that $\mathfrak{S}^g$ separates $N_\smallinfinity^g$ into two components $A$ and $E^g$: the first the geometric limit of the $A_n$; the second the geometric limit of the $E_n$.

\begin{lem}
The subgroup $\Gamma_\smallinfinity^a$ carries the fundamental group of $E^g$, and so $E^g$ is an end of $N_\smallinfinity^g$ that lifts isometrically to $N_\smallinfinity^a$.
\end{lem}

We establish some notation before proving the lemma.

Let $\mathcal{C}^a$ be the component of the boundary of the convex core of $N_\smallinfinity^a$ facing the end $E^a$, and let $\mathcal{C}^g$ be its image in $N_\smallinfinity^g$.
The surfaces $\mathcal{C}^g$ and $\mathfrak{S}^g$ are homotopic in $N_\smallinfinity^g$, and we fix a homotopy.

We consider the component of the convex core boundary of $N_n$ facing $\sigma_M(M_n)$ as a pleated surface $\mathcal{C}_n \to N_n$. 
When no confusion can arise, we blur the distinction between $\mathcal{C}_n \to N_n$ and its image.
Each $\mathcal{C}_n$ is compact since $\sigma_M(M_n)$ lies in $\T(\overline{\partial_0 M})$.

After identifying the $\K_n$ with $M$, the manifolds $N_n$ are each marked by the unique \textit{bona fide} lift $\partial_0 M \to N_n$, which, in turn, marks each $\mathcal{C}_n$.

The pleated surfaces $\mathcal{C}_n \to N_n$ descend to pleated surfaces $\mathcal{D}_n \to M_n$.

\begin{proof}[Proof of lemma]
We show that every essential closed curve in $N_\smallinfinity^g$ is homotopic into $A$.

\medskip \noindent
\textsc{Case I.} Suppose that, after passing to a subsequence, there is a compact set $K$ in $N_\smallinfinity^g$ whose preimage $K_n$ in the approximate $N_n$ intersects $\mathcal{C}_n$ for all $n$.
By Theorem \ref{pleatedcompact} (Pleated Surfaces Compact) and geometric convergence of the $N_n$,  we see, after choosing base frames and taking a subsequence, that the $\mathcal{C}_n \to N_n$ are converging to a pleated surface $\mathcal{C} \to N_\smallinfinity^g$.

If $\mathcal{C}$ is noncompact, then there is a fixed curve $\gamma$ in $\partial_0 M$ whose length in $\mathcal{C}_n$ is tending to zero.
This implies that the end $E^a$ has a rank--one cusp, which we have excluded. 
To see this, let $\gamma_n$ be the geodesic representative of $\gamma$ in $\mathcal{C}_n$, and let $\gamma_n^*$ be its geodesic representative in $N_n$.
The cores $\K_n$ have uniformly bounded diameters, so we may choose an $\epsilon$ so that the $\epsilon$--Margulis tubes $\mathbf{T}_n$ about the $\gamma_n^*$ miss the $\K_n$, by Brooks and Matelski's theorem \cite{brooksmatelski}.
Since the lengths of the $\gamma_n$ are tending to zero, the depths at which they lie in the $\mathbf{T}_n$ must be tending to infinity, again by \cite{brooksmatelski}.
It follows that $\partial \mathbf{T}_n$ intersects the boundary of a large neighborhood of $\mathcal{C}_n$ in the end of $N_n$ facing $\sigma_M(M_n)$.
Since $\K_n$ misses $\mathbf{T}_n$, the distance from $\K_n$ to $\mathcal{C}_n$ is uniformly bounded, and the $\K_n$ have uniformly bounded diameters, the $\K_n$ must eventually lie to the left of the $\mathbf{T}_n$.
So the $\gamma_n^*$ eventually lie to the right of $\K_n$, and we discover a rank--one cusp in $E^a$.

So $\mathcal{C}$ is compact.
It follows that the $\mathcal{C}_n$ have uniformly bounded diameter, and so, for sufficiently large $n$, we may push the $\mathcal{C}_n$ into $N_\smallinfinity^g$ to obtain surfaces $\mathcal{C}_n^g$ converging to $\mathcal{C}$.

Moreover, since the $\Gamma_n$ converge algebraically, the surface $\mathcal{C}$ is in the same homotopy class as $\mathcal{C}^g$, and we fix a homotopy between them.
Since the $\mathcal{C}_n^g$ are converging to $\mathcal{C}$,  the surfaces $\mathcal{C}_n^g$ and $\mathcal{C}$ are eventually homotopic in the $1$--neighborhood of $\mathcal{C}$.

It follows that the $\mathcal{C}_n$ admit homotopies to the $\mathfrak{S}_n$ of uniformly bounded diameter: the $\mathcal{C}_n^g$ are uniformly homotopic to $\mathcal{C}$ in $N_\smallinfinity^g$, which is homotopic to $\mathcal{C}^g$, which is homotopic to $\mathfrak{S}^g$; and we may push the resulting homotopy from $\mathcal{C}_n^g$ to $\mathfrak{S}^g$ back into the approximates.

Now, let $\gamma$ be a closed geodesic in $N_\smallinfinity^g$, and push it back to curves $\gamma_n$ in the approximates.
The geometric convergence implies that the geodesic curvatures of the $\gamma_n$ are all uniformly very nearly zero in the $N_n$ for all large $n$.
It follows that each $\gamma_n$ is homotopic to its geodesic representative $\gamma_n^*$ in the $\epsilon$--neighborhood of $\gamma_n^*$.
Now, each $\gamma_n^*$ lies in the convex core of $N_n$, to the left of the pleated surface $\mathcal{C}_n$.
Since the $\mathcal{C}_n$ admit homotopies to the $\mathfrak{S}_n$ of uniformly bounded diameter, there are ambient homotopies in the $N_n$ with supports of uniformly bounded diameter that carry the $\gamma_n^*$ to the left of $\mathfrak{S}_n$.
Since the $\epsilon$--neighborhoods of the $\gamma_n^*$ eventually map to the geometric limit, we conclude that $\gamma$ is homotopic into $A$.

So any closed geodesic in $N_\smallinfinity^g$ is homotopic into $A$.

Now let  $\mathbf P$ be a cusp in $N_\smallinfinity^g$,
let $p$ be a parabolic element in the corresponding subgroup of $\Gamma_\smallinfinity^g$,
and let $\eta$ be a geodesic representing $p$ in the flat metric on $\partial \mathbf P$.

Suppose that there is a sequence $\{p_n \, | \, p_n \in \Gamma_n\}$ of loxodromic elements  converging to $p$.
Then the cusp $\mathbf P$ is the geometric limit of Margulis tubes $\mathbf T_n$ about the corresponding geodesics $\gamma_n^*$ in $N_n$.
Pushing $\eta$ into the approximates yields a sequence of closed curves $\eta_n$ whose geodesic curvatures are close to one, and which are hence uniformly close to the geodesic $\eta_n^*$ in $\partial \mathbf T_n$ in the same homotopy class.

The depths of the tubes $\mathbf T_n$ are tending to infinity.
The geodesic $\gamma_n^*$ lies in the convex core of $N_n$, to the left of $\mathcal{C}_n$, and since the $\mathcal{C}_n$ have uniformly bounded diameter, and hence cannot penetrate too deeply into $\mathbf T_n$, there is a smaller Margulis tube about $\gamma_n^*$ that lies to the left of $\mathcal{C}_n$ and whose boundary is a uniformly bounded distance from $\eta_n$.

If $p$ is the limit of parabolic elements $p_n$, then $\mathbf P$ is the geometric limit of Margulis tubes $\mathbf T_n$, and, as above, we may take these to lie in the convex cores of the approximates, perhaps after changing $\mathbf P$ slightly.

In either case, we obtain bounded diameter homotopies carrying the $\eta_n$ to the left of the $\mathcal{C}_n$, and then ambient homotopies with supports of uniformly bounded diameter carrying the result to the left of $\mathfrak{S}_n$.
We conclude that $\eta$ is homotopic into $A$.

So, in \textsc{Case I}, every essential closed curve in $N_\smallinfinity^g$ is homotopic into $A$.

\medskip \noindent
\textsc{Case II.} Suppose that for every compact set $K$ in $N_\smallinfinity^g$, there are only finitely many $n$ such that the preimage $K_n$ of $K$ in $N_n$ intersects $\mathcal{C}_n$.

It follows that the images $\mathcal{D}_n$ of the $\mathcal{C}_n$ in $M_n$ must eventually lie outside of the compact cores $\K_n$.

To see this, suppose to the contrary that, after passing to a subsequence, each $\mathcal{D}_n$ intersects $\K_n$.    
Then, after choosing base frames and passing to another subsequence, the $\mathcal{D}_n \to M_n$ converge to a finite area pleated surface $\mathcal{D} \to M_\smallinfinity^g$ freely homotopic into $\partial \K^g$ in $M_\smallinfinity^g$.
So, if $\mathcal{G}$ is a closed geodesic in $\mathcal{D}$, 
and $\mathcal{G}_n$ a sequence of geodesics in $\mathcal{D}_n$ converging to $\mathcal{G}$, 
then there is a uniform $R>0$ such that the $\mathcal{G}_n \to M_n$ are eventually homotopic into $\partial \K_n$ in the $R$--neighborhood of the latter.

Now, since $\mathcal{D}_n$ is peripheral and $M$ is acylindrical, the surface $\mathcal{C}_n \to N_n$ is the unique lift of $\mathcal{D}_n \to M_n$ that is not simply connected.                                                                
So, lifting the homotopies of $\mathcal{G}_n \to M_n$ into $\partial \K_n$ to $N_n$, we find that each $\mathcal{C}_n$ intersects the $R$--neighborhood of $\mathfrak{S}_n$.
But the $\mathfrak{S}_n$ are the preimages of the compact set $\mathfrak{S}^g$ in $N_\smallinfinity^g$.

We conclude that the $\mathcal{D}_n \to M_n$ eventually miss the cores $\K_n$.
Since they cannot lie in the rank--two cusps of the $M_n$, they eventually lie to the left of the $\partial \K_n$.
It follows that the $\mathcal{C}_n$ eventually lie to the left of the $\mathfrak{S}_n$.
(In the case where $\partial_0 M$ is disconnected, the $\mathcal{D}_n$ are all homotopic to a fixed component of $\partial_0 M$, and as distinct components of $\partial_0 M$ are not homotopic, one concludes that the $\mathcal{D}_n$ lie to the left of the chosen one.)

Again, let $\gamma$ be a closed geodesic in $N_\smallinfinity^g$, and push it back to curves $\gamma_n$ in the $N_n$.
As before, the $\gamma_n$ are homotopic to their geodesic representatives $\gamma_n^*$ in the $\epsilon$--neighborhood of $\gamma_n^*$.
But now, $\gamma_n^*$ lies in the convex core of $N_n$ which lies to the left of $\mathfrak{S}_n$, and, going back to the limit, we conclude that $\gamma$ is homotopic into $A$.

Similarly, the argument that the parabolic elements are homotopic into $A$ proceeds as in  \textsc{Case I}, the argument simplified by $\mathcal{C}_n$ lying to the left of $\mathfrak{S}_n$.

\medskip \noindent
We have now shown that any essential closed curve in $N_\smallinfinity^g$ is homotopic into $A$.
It follows that $\pi_1(\mathfrak{S}^g) \cong \Gamma_\smallinfinity^a < \Gamma_\smallinfinity^g$ carries the fundamental group of $E^g$, and so $E^g$ is an end of $N_\smallinfinity^g$ that lifts isometrically to $E^a$.
\end{proof}

To complete the proof, choose an Epstein surface $\F$ in $E^g$ (see Section \ref{usefulsurfaces}), and push it into the approximates to obtain surfaces $\F_n$. 
By the geometric convergence, for all large $n$ the surfaces $\F_n$ are strictly convex and the normal curvatures of the $\F_n$ converge to those of $\F$.  
Paired with the isometry $E^g \to E^a$, convexity provides normal projections
\[
\sigma_M(M_n) \to \F_n 
\quad \mathrm{and} \quad
\sigma_M(M_\smallinfinity^a) \to \F
\]
whose derivatives depend only on the normal curvatures of the $\F_n$ and $\F$. 
Composing with the approximating bilipschitz maps, we see that the derivatives of the normal projections will ``cancel" in the limit, as the limiting bilipschitz map is the restriction of an isometry (which will carry principal directions to principal directions)---as we are only concerned with continuity here, we do pass to the limit (compare the proofs of Theorems \ref{collar} and \ref{fill}).
By $C^{\, \smallinfinity} \! \!$--convergence, we conclude that for large $n$ the composition is very close to conformal and so the $\sigma_M(M_n)$ are converging to $\sigma_M(M_\smallinfinity^a)$.
\end{proof}


\bibliographystyle{plain}
\bibliography{bartholomew}

\medskip

\noindent Department of Mathematics, Brown University, Providence, RI 02912 \newline \noindent  \texttt{rkent@math.brown.edu}

\end{document}